arXiv:1507.08323v1 [math.FA] 29 Jul 2015

# CONVEX-CYCLIC MATRICES, CONVEX-POLYNOMIAL INTERPOLATION & INVARIANT CONVEX SETS

NATHAN S. FELDMAN & PAUL MCGUIRE

Abstract. We define a convex-polynomial to be one that is a convex combination of the monomials $\{1, z, z^2, \ldots\}$. This paper explores the intimate connection between peaking convex-polynomials, interpolating convex-polynomials, invariant convex sets, and the dynamics of matrices. In particular, we use these intertwined relations to both prove which matrices are convex-cyclic while at the same time proving that we can prescribe the values and a finite number of the derivatives of a convex-polynomial subject to certain natural constraints. These properties are also equivalent to determining those matrices whose invariant closed convex sets are all invariant subspaces.

Our characterization of the convex-cyclic matrices gives a new and correct proof of a similar result by Rezaei that was stated and proven incorrectly.

## 1. Introduction

If $E$ is a subset of a vector space $X$, then the convex hull of $E$, denoted by $co(E)$, is the set of all convex combinations of elements of $E$; that is, all finite linear combinations of the elements of $E$ where the coefficients are non-negative and sum to one. We will let $\mathcal{CP}$ denote the convex hull of the set of monomials $\{1, z, z^2, z^3, \ldots\}$ within the vector space of all polynomials in $z$. Thus,

$$\mathcal{CP} = co(\{1, z, z^2, z^3, \ldots\}) = \left\{ \sum_{k=0}^{n} a_k z^k : a_k \geq 0 \text{ for all } 0 \leq k \leq n \text{ and } \sum_{k=0}^{n} a_k = 1 \right\}.$$

Following Rezaei [12], we will call the elements of $\mathcal{CP}$ *convex-polynomials.*

If $T$ is a continuous linear operator on a locally convex-space $X$ and $x \in X$, then the orbit of $x$ under $T$ is $Orb(x, T) = \{T^n x\}_{n=0}^{\infty} = \{x, Tx, T^2x, \ldots\}$. A continuous linear operator $T$ on $X$ is said to be *cyclic* if there is a vector $x \in X$ such that the *linear span* of the orbit of $x$ under $T$ is dense in $X$; that is if $\{p(T)x : p \text{ is a polynomial}\}$ is dense in $X$. Following Rezaei [12], we define an operator $T$ to be *convex-cyclic* if there is a vector $x \in X$ such that the *convex-hull* of the orbit of $x$ under $T$ is dense in $X$; that is if $\{p(T)x : p \in \mathcal{CP}\}$ is dense in $X$. Convex-cyclic operators were introduced by Rezaei [12] and have been studied in [3] and[11]. More generally the dynamics of matrices have been studied in [1], [5], [6] and in their references.

While the problem of determining which matrices are convex-cyclic was addressed in Rezaei [12], an unfortunate oversight was made resulting in an incorrect statement and proof of the result. While we deduce similar results without the oversight, our focus remains on the development of a framework for working with convex polynomials and matrices and, as a consequence, our proofs and approach are markedly different from those in Rezaei [12]. Furthermore our techniques can be applied in infinite dimensional settings as well.

We will prove three main theorems. The first theorem describes which matrices are convex-cyclic. The characterization is natural and simple: A matrix is convex-cyclic exactly when it is cyclic and satisfies some necessary conditions on its eigenvalues. Recall that a matrix is cyclic if and only if each of its eigenvalues has geometric multiplicity one.

The second theorem says that under natural necessary conditions we can prescribe the values and a finite number of derivatives of a convex-polynomial.

The third theorem gives a simple condition on the eigenvalues of a matrix so that each of its invariant closed convex sets is actually an invariant subspace. This condition has been given by Elsner [7]. However we use fundamentally different techniques in our approach, that may also be applied in infinite dimensions.

**Theorem 1.1.** ***Convex-Cyclicity of Matrices***

***The Real Case:*** *If $T$ is a real $n \times n$ matrix, then $T$ is convex-cyclic on $\mathbb{R}^n$ if and only if $T$ is cyclic and its real and complex eigenvalues are contained in $\mathbb{C} \setminus (\overline{\mathbb{D}} \cup \mathbb{R}^+)$. If $T$ is convex-cyclic, then the convex-cyclic vectors for $T$ are the same as the cyclic vectors for $T$ and they form a dense set in $\mathbb{R}^n$.*

***The Complex Case:*** *If $T$ is an $n \times n$ matrix, then $T$ is convex-cyclic on $\mathbb{C}^n$ if and only if $T$ is cyclic and its eigenvalues $\{\lambda_k\}_{k=1}^n$ are all contained in $\mathbb{C} \setminus (\overline{\mathbb{D}} \cup \mathbb{R})$ and satisfy $\lambda_j \neq \overline{\lambda}_k$ for all $1 \leq j, k \leq n$. If $T$ is convex-cyclic, then the convex-cyclic vectors for $T$ are the same as the cyclic vectors for $T$ and they form a dense set in $\mathbb{C}^n$.*

**Theorem 1.2.** ***Convex-Polynomial Interpolation***

*If $S = \{x_k\}_{k=1}^m \cup \{z_k\}_{k=1}^n \subseteq \mathbb{C}$ where $\{x_k\}_{k=1}^m \subseteq \mathbb{R}$ and $\{z_k\}_{k=1}^n \subseteq \mathbb{C} \setminus \mathbb{R}$, then the following are equivalent:*

*(a) for any finite set $\{y_{j,k} : 0 \leq j \leq N, 1 \leq k \leq m\}$ of real numbers and for any finite set $\{w_{j,k} : 0 \leq j \leq N, 1 \leq k \leq n\}$ of complex numbers there exists a convex-polynomial $p$ such that $p^{(j)}(x_k) = y_{j,k}$ for all $0 \leq j \leq N$ and $1 \leq k \leq m$ and $p^{(j)}(z_k) = w_{j,k}$ for all $0 \leq j \leq N$ and $1 \leq k \leq n$.*

*(b) The real numbers $\{x_k\}_{k=1}^m$ are distinct and satisfy $\{x_k\}_{k=1}^m \subseteq (-\infty, -1)$ and the imaginary numbers $\{z_k\}_{k=1}^n$ are distinct, $\{z_k\}_{k=1}^n \subseteq \mathbb{C} \setminus \overline{\mathbb{D}}$ and $z_j \neq \overline{z}_k$ for all $1 \leq j, k \leq n$.*

**Theorem 1.3.** ***Invariant Convex sets for Matrices***

***The Complex Case:*** *A matrix $T$ acting on $\mathbb{C}^n$ has the property that all of its invariant closed convex-sets are invariant subspaces if and only if the eigenvalues $\{\lambda_k\}_{k=1}^n$ of $T$ are all contained in $\mathbb{C} \setminus (\overline{\mathbb{D}} \cup \mathbb{R})$ and satisfy $\lambda_j \neq \overline{\lambda}_k$ for all $1 \leq j, k \leq n$.*

***The Real Case:*** *A matrix $T$ acting on $\mathbb{R}^n$ has the property that all of its invariant closed convex-sets are invariant subspaces if and only if all of its (real and complex) eigenvalues are contained in $\mathbb{C} \setminus (\overline{\mathbb{D}} \cup \mathbb{R}^+)$.*

These three theorems above are actually all equivalent and we will begin by proving the first two of them simultaneously in a series of intertwined steps. First we will construct some convex-polynomials which "peak" at a prescribed point in a given finite set of points. This will then allow us to easily determine which diagonal matrices are convex-cyclic, from which we will be able to prove an interpolation theorem for the values (but no derivatives) of a convex-polynomial. That interpolation theorem will then be used to prove that certain non-diagonalizable matrices are convex-cyclic, which will give a stronger interpolation theorem. This will, in turn, lead to a larger class of convex-cyclic matrices, and so forth. This continues until we arrive at the above three theorems.

## 2. Convex-Polynomials & Necessary Conditions

Let $f^{(k)}$ denote the $k^{th}$ derivative of a function $f$. If $f$ has a power series that converges near zero, say $f(z) = \sum_{k=0}^{\infty} a_k z^k$, then the coefficients of the series are related to the derivatives of $f$ by $a_k = f^{(k)}(0)/k!$ for $k \geq 0$. This fact can be used to easily characterize the convex-polynomials and determine some basic properties of them, as stated in the following proposition. Recall that $\mathcal{CP}$ denotes the set of convex-polynomials.

**Proposition 2.1.** *The following hold:*

1. *A polynomial $p(z)$ is a convex-polynomial if and only if $p^{(k)}(0) \geq 0$ for all $k \geq 0$ and $p(1) = 1$.*
2. *The set $\mathcal{CP}$ is closed under composition and multiplication.*
3. *If $p$ is a convex-polynomial, then $p(\mathbb{R}) \subseteq \mathbb{R}$, $p(\overline{\mathbb{D}}) \subseteq \overline{\mathbb{D}}$, and $\overline{p(z)} = p(\overline{z})$ for all $z \in \mathbb{C}$. In particular, $|p(z)| \leq 1$ whenever $|z| \leq 1$ and $p \in \mathcal{CP}$.*

2.1. **Necessary Conditions.** Before proceeding, it is useful to provide three simple but illustrative examples highlighting some necessary conditions.

**Example 2.2.** Let $T = \begin{bmatrix} \lambda_1 & 0 \\ 0 & \lambda_2 \end{bmatrix}$ be a matrix on $\mathbb{F}^2$ where $\mathbb{F}$ equals $\mathbb{R}$ or $\mathbb{C}$.

1. If $T$ is convex-cyclic, then $\lambda_1 \neq \lambda_2$.
2. If $T$ is convex-cyclic, then $|\lambda_i| > 1$ for $i = 1, 2$.
3. If $T$ is convex-cyclic, then $\lambda_2 \neq \overline{\lambda_1}$.

*Proof.* If $v = \begin{bmatrix} v_1 \\ v_2 \end{bmatrix}$ and $p$ is a polynomial, then $p(T)v = \begin{bmatrix} p(\lambda_1)v_1 \\ p(\lambda_2)v_2 \end{bmatrix}$.

(1) If $\lambda_1 = \lambda_2 = \lambda$, then for any vector $v$, we have $p(T)v = p(\lambda)v$. Thus $\{p(T)v : p \text{ is a polynomial}\}$ is contained in a one-dimensional subspace and is therefore not dense in $\mathbb{F}^2$. It follows that $T$ is not cyclic, so $T$ is not convex-cyclic either.

(2) If $|\lambda_i| \leq 1$ for some $i$, then by Proposition 2.1 for any $p \in \mathcal{CP}$ we have $|p(\lambda_i)| \leq 1$ so $|p(\lambda_i)v_i| \leq |v_i|$. Thus, $\{p(T)v : p \in \mathcal{CP}\}$ is not dense in $\mathbb{F}^2$ for any vector $\vec{v}$. Thus $T$ is not convex-cyclic.

(3) If $\lambda_2 = \overline{\lambda_1}$ and $p \in \mathcal{CP}$, then by Proposition 2.1 we have $p(\lambda_2) = p(\overline{\lambda_1}) = \overline{p(\lambda_1)}$. Thus, $p(T)v = \begin{bmatrix} p(\lambda_1)v_1 \\ p(\lambda_2)v_2 \end{bmatrix} = \begin{bmatrix} p(\lambda_1)v_1 \\ \overline{p(\lambda_1)}v_2 \end{bmatrix}$. In particular, the vector $\begin{bmatrix} 1 \\ 0 \end{bmatrix}$ is not in the closure of $\{p(T)v : p \in \mathcal{CP}\}$. So, $T$ is not convex cyclic. □

We will show, in Theorem 5.1, that excluding the above three simple situations, all other diagonal matrices are convex cyclic. In the paper by H. Rezaei [12], the third condition above was overlooked.

We next give some necessary conditions that apply not only to matrices, but to operators on locally convex spaces. Items (2) and (3) below in Proposition 2.3 are more general versions of items (2) and (3) from Example 2.2 above.

Let $\sigma_p(T)$ denote the point spectrum of the operator $T$, that is, the set of all eigenvalues of $T$.

**Proposition 2.3.** ***Necessary Conditions on Eigenvalues of** $T^*$.* *If $T$ is a convex-cyclic continuous linear operator on a complex locally convex vector space $X$ over $\mathbb{C}$, then the following hold:*

1. *$\sigma_p(T^*) \cap \mathbb{R} = \emptyset$.*
2. *$\sigma_p(T^*) \cap \overline{\mathbb{D}} = \emptyset$.*
3. *If $\lambda_1, \lambda_2 \in \sigma_p(T^*)$, then $\lambda_2 \neq \overline{\lambda_1}$.*

*Proof.* Let $x$ be a convex-cyclic vector for $T$, thus $\{p(T)x : p \in \mathcal{CP}\}$ is dense in $X$. Also let $\lambda$ be an eigenvalue for $T^*$ with eigenvector $v \in X^*$. So $T^*v = \lambda v$. It follows that for any polynomial $p$, $p(T^*)v = p(\lambda)v$. If $p$ is a convex-polynomial, then $p$ has real coefficients so $p(T)^* = p(T^*)$. So when $p \in \mathcal{CP}$ we have

$$(*) \qquad \langle p(T)x, v\rangle = \langle x, p(T)^*v\rangle = \langle x, p(T^*)v\rangle = \langle x, p(\lambda)v\rangle = \overline{p(\lambda)} \cdot \langle x, v\rangle.$$

Since $\{p(T)x : p \in \mathcal{CP}\}$ is dense in $X$, it follows that $\{\overline{p(\lambda)}\langle x, v\rangle : p \in \mathcal{CP}\}$ must be dense in $\mathbb{C}$. However, for any convex-polynomial $p$ we have that $p(\mathbb{R}) \subseteq \mathbb{R}$ and $p(\overline{\mathbb{D}}) \subseteq \overline{\mathbb{D}}$, so it follows that if $\lambda \in \mathbb{R} \cup \overline{\mathbb{D}}$, then $\{\overline{p(\lambda)}\langle x, v\rangle : p \in \mathcal{CP}\}$ is not dense in $\mathbb{C}$. Thus we (1) and (2) hold.

To see that (3) holds suppose, by way of contradiction, that $\lambda_1, \lambda_2 \in \sigma_p(T^*)$ and that $\lambda_2 = \overline{\lambda_1}$. Let $\lambda := \lambda_2 = \overline{\lambda_1}$. Then both $\lambda$ and $\overline{\lambda}$ are eigenvalues for $T^*$. Let $v_1$ and $v_2$ be eigenvectors for $T^*$ with eigenvalues $\lambda$ and $\overline{\lambda}$ respectively. By (1) we know that $\lambda \notin \mathbb{R}$, thus $\lambda \neq \overline{\lambda}$. It follows that $v_1$ and $v_2$ are linearly independent. From $(*)$ above we know that for every $p \in \mathcal{CP}$ we have

$$(**) \qquad \begin{bmatrix} \langle p(T)x, v_1\rangle \\ \langle p(T)x, v_2\rangle \end{bmatrix} = \begin{bmatrix} \overline{p(\lambda)} \cdot \langle x, v_1\rangle \\ \overline{p(\overline{\lambda})} \cdot \langle x, v_2\rangle \end{bmatrix} = \begin{bmatrix} \overline{p(\lambda)} \cdot \langle x, v_1\rangle \\ p(\lambda) \cdot \langle x, v_2\rangle \end{bmatrix}$$

where $\overline{p(\overline{\lambda})} = p(\lambda)$ since $p$ has real coefficients. Now since $\{p(T)x : p \in \mathcal{CP}\}$ is dense in $X$ and $v_1$ and $v_2$ are linearly independent it follows that as $p$ varies over all convex-polynomials that the left hand side of $(**)$ is dense in $\mathbb{C}^2$. However, as $p$ varies over all convex-polynomials the right hand side of $(**)$ is not dense in $\mathbb{C}^2$ (if one coordinate is small in absolute value the other coordinate will also be small). This is a contradiction. It follows that (3) holds. □

## 3. Dense Convex Sets & The Hahn-Banach Theorem

The following result is a fundamental tool for studying convex-cyclic operators.

**Theorem 3.1** (A Criterion for a Convex Set to be Dense)**.** *If $C$ is a convex set in a locally convex linear space $X$, then $C$ is dense in $X$ if and only if for every non-zero continuous linear functional $f$ on $X$ we have that $\sup_{x\in C} Re(f(x)) = \infty$. Furthermore if $S \subseteq C$ and the convex-hull of $S$ is dense in $C$, then*

$$\sup_{x\in C} Re(f(x)) = \sup_{x\in S} Re(f(x)).$$

The previous result is a simple consequence of the geometric form of the Hahn-Banach Theorem which says that whenever a point does not belong to a closed convex set, then the point and the convex set can be strictly separated by a real hyperplane. See [4, Theorem 3.13, p. 111].

A vector $x \in X$ is a *convex-cyclic vector* for $T$ if $co(Orb(x, T))$ is dense in $X$.

**Corollary 3.2** (Hahn-Banach Characterization of Convex-Cyclicity)**.** *Let $X$ be a locally convex space over the real or complex numbers, $T : X \to X$ a continuous linear operator, and $x \in X$. Then the following are equivalent:*

(1) *The vector $x$ is a convex-cyclic vector for $T$.*
(2) *For every non-zero continuous linear functional $f$ on $X$ we have*

$$\sup_{n\geq 0} Re[f(T^n x)] = \infty.$$

(3) *For every non-zero continuous linear functional $f$ on $X$ we have*

$$\sup\{Re[f(p(T)x)] : p \in \mathcal{CP}\} = \infty.$$

*Proof.* Apply Theorem 3.1 where $C$ is the convex hull of the orbit of $x$. □

Next we use the above Hahn-Banach characterization from above to establish a condition for the direct sum of two convex-cyclic operators to be convex-cyclic.

A set $E$ is *bounded* in a locally convex-space $X$ if for every neighborhood $U$ of zero there is a $c > 0$ such that $E \subseteq cU$. This is equivalent to being weakly-bounded which says that $f(E)$ is a bounded set of scalars for every continuous linear functional $f$ on $X$.

We will say that a continuous linear operator $T$ on a locally convex-space is *power bounded* if all of the orbits of $T$ are bounded sets. This is consistent with the notion of power boundedness on a Banach space.

Let $\mathbb{F}$ denote either $\mathbb{R}$ or $\mathbb{C}$.

**Proposition 3.3.** ***Direct Sums of Convex-Cyclic Operators.***

*Let $T_1$ and $T_2$ be continuous convex-cyclic linear operators on locally convex spaces $X_1$ and $X_2$ over $\mathbb{F}$. If there exists a convex-polynomial $p_0$ such that $p_0(T_1)$ is convex-cyclic and $p_0(T_2)$ is power bounded, then $T_1 \oplus T_2$ is convex-cyclic on $X_1 \oplus X_2$.*

*Furthermore, if $u_1$ is a convex-cyclic vector for $p_0(T_1)$ and $u_2$ is a convex-cyclic vector for $T_2$, then $\vec{u} = (u_1, u_2) \in X_1 \oplus X_2$ is a convex-cyclic vector for $T_1 \oplus T_2$.*

Notation: If $x \in X$ and $f \in X^*$ we will use both $\langle x, f\rangle$ and $f(x)$ to denote the value of $f$ acting on the vector $x$.

*Proof.* Let $T = T_1 \oplus T_2$, $X = X_1 \oplus X_2$, and let $p_0$ be a convex-polynomial such that $p_0(T_1)$ is convex-cyclic and $p_0(T_2)$ is power bounded. Also let $\vec{f} = (f_1, f_2) \in X^* \setminus \{0\} = (X_1^* \oplus X_2^*) \setminus \{0\}$ and let $\vec{u} = (u_1, u_2) \in X$ be a vector where $u_1$ is a convex-cyclic vector for $p_0(T_1)$ and $u_2$ is a convex-cyclic vector for $T_2$. We must show that $\sup\limits_{p \in \mathcal{CP}} Re\langle p(T)\vec{u}, \vec{f}\rangle = \infty$.

*Case 1:* $f_1 = 0$.

In this case, note that $f_2 \neq 0$ since $\vec{f} \neq \vec{0}$, and thus we have

$$\sup_{p \in \mathcal{CP}} Re\langle p(T)\vec{u}, \vec{f}\rangle = \sup_{p \in \mathcal{CP}} Re\left[\langle p(T_1)u_1, f_1\rangle + \langle p(T_2)u_2, f_2\rangle\right] = \sup_{p \in \mathcal{CP}} Re\langle p(T_2)u_2, f_2\rangle = \infty. \tag{1}$$

The last supremum above is infinite by Corollary 3.2, since $T_2$ is convex-cyclic, $u_2$ is a convex-cyclic vector for $T_2$ and $f_2 \neq 0$.

*Case 2:* $f_1 \neq 0$.

Again we must show that $\sup\limits_{p \in \mathcal{CP}} Re\langle p(T)\vec{u}, \vec{f}\rangle = \infty$. Using the fact that $T_1$ is convex-cyclic and $p_0(T_2)$ is power bounded we have

$$\sup_{p \in \mathcal{CP}} Re\langle p(T)\vec{u}, \vec{f}\rangle \geq \sup_{n \geq 1} Re\langle p_0(T)^n\vec{u}, \vec{f}\rangle =$$

$$\sup_{n \geq 1}[Re\langle p_0(T_1)^n u_1, f_1\rangle + Re\langle p_0(T_2)^n u_2, f_2\rangle] = \infty.$$

The last equality above holds by Corollary 3.2 since $p_0(T_1)$ is convex-cyclic and $p_0(T_2)$ is power bounded. The last condition implies that $\sup_{n \geq 1} |\langle p_0(T_2)^n x, \vec{f}\rangle| < \infty$. Also since $u_1$ is a convex-cyclic vector for $p_0(T_1)$ and $f_1 \neq 0$ then by Corollary 3.2 we have that $\sup_{n \geq 1} Re\langle p_0(T_1)^n u_1, f_1\rangle = \infty$. The theorem now follows. □

The previous theorem will be useful in showing that block diagonal matrices are convex-cyclic when each block is convex-cyclic.

## 4. Peaking Convex-Polynomials

If $T \subseteq \mathbb{C}$ and $f : T \to \mathbb{C}$ is a bounded function defined on $T$, then we say that $f$ *peaks* on $T$ if $|f|$ attains it supremum on $T$ at a unique point in $T$. That is, $f$ peaks on $T$ if there exists an $x_0 \in T$ such that $|f(x_0)| > |f(x)|$ for all $x \in T \setminus \{x_0\}$. If $S \subseteq T$, then we will say that $f : T \to \mathbb{C}$ peaks on $T$ at a point in $S$ if there exists an $x_0 \in S$ such that $|f(x_0)| > |f(x)|$ for all $x \in T \setminus \{x_0\}$.

For a bounded function $f : T \to \mathbb{C}$ let $\|f\|_T := \sup\{|f(x)| : x \in T\}$.

For $0 \le \alpha \le 1$, let

$$p_\alpha(z) = \alpha z + (1 - \alpha)$$

and notice that $p_\alpha$ is a convex polynomial. Also, given a non-negative integer $m$ and $0 \le \alpha \le 1$ define the polynomials

$$p_{m,\alpha}(z) = z^m p_\alpha(z) = z^m(\alpha z + (1 - \alpha)) = \alpha z^{m+1} + (1 - \alpha)z^m.$$

Note that $p_{m,\alpha}$ is a convex-polynomial since it is a product of two convex-polynomials, also because it is simply a convex combination of $z^m$ and $z^{m+1}$.

**Theorem 4.1.** ***Peaking Convex-Polynomials.***
*If $S = \{z_k\}_{k=1}^N$ is a finite set of complex numbers satisfying:*

1. *the points in $S$ are distinct;*
2. *$R = \max\{|z| : z \in S\} > 1$; and*
3. *$z_j \ne \overline{z}_k$ whenever $|z_j| = |z_k| = R$ and $j \ne k$,*

*then the following holds:*

*for every $\alpha \in (0, 1)$, except possibly one, there exists a $K > 0$ such that for all $m \ge K$ the convex-polynomial $p_{m,\alpha}(z) = z^m(\alpha z + 1 - \alpha)$ peaks on $S$ at some point $z_{k_0} \in S$ satisfying $|z_{k_0}| = R$ and $\|p_{m,\alpha}\|_S = c_\alpha R^m$ where $c_\alpha = \max\{|\alpha z_k + (1 - \alpha)| : |z_k| = R\}$. Moreover if $S \cap \mathbb{R} = \emptyset$ and $m \ge 1$, then for all but finitely many $\alpha \in (0, 1)$ we have $p_{m,\alpha}(S) \cap \mathbb{R} = \emptyset$.*

*Proof.* Let $S = \{z_k\}_{k=1}^N$, $R = \|z\|_S = \max\{|z| : z \in S\}$, $T_1 = \{z_k : |z_k| = R\}$, and $T_2 = S \setminus T_1 = \{z_k : |z_k| < R\}$. Note that by (2) in our hypothesis we have that $R > 1$.

Let $0 < \alpha < 1$ and $p_\alpha(z) = \alpha z + (1 - \alpha)$ and note that $|p_\alpha(z)| = \alpha \cdot \left|z - \frac{\alpha-1}{\alpha}\right|$ is $\alpha$ times the distance from $z$ to $\frac{\alpha-1}{\alpha}$ and $\frac{\alpha-1}{\alpha} < 0$. Consider $\|p_\alpha\|_{T_1}$ and $\|p_\alpha\|_{T_2}$. In both cases $|p_\alpha|$ attains its absolute maximum on the given set at a point of that set that is the greatest distance away from $\frac{\alpha-1}{\alpha}$. If $z_k \in T_1$ and $z_k = Re^{i\theta_k}$, then

$$|p_\alpha(z_k)|^2 = |\alpha z_k + (1 - \alpha)|^2 = |\alpha Re^{i\theta_k} + 1 - \alpha|^2 =$$

$$\Big(\alpha R \cos(\theta_k) + 1 - \alpha\Big)^2 + \alpha^2 R^2 \sin^2(\theta_k) =$$

$$\alpha^2 R^2 \cos^2(\theta_k) + 2\alpha R \cos(\theta_k)(1 - \alpha) + (1 - \alpha)^2 + \alpha^2 R^2 \sin^2(\theta_k) =$$

$$\alpha^2 R^2 + (1 - \alpha)^2 + 2\alpha(1 - \alpha)R \cos(\theta_k).$$

Recall $z_j \ne \overline{z}_k$ for all $1 \le j, k \le N$ with $j \ne k$ and $|z_j| = |z_k| = R$. Thus, for any $z_j, z_k \in T_1$ with $j \ne k$ we must have $Re(z_j) \ne Re(z_k)$, thus $R\cos(\theta_k) = Re(z_k) \ne Re(z_j) = R\cos(\theta_j)$. Hence there is a unique $z_{k_0} \in T_1$, independent of $\alpha$, for which $2\alpha(1 - \alpha)R\cos(\theta_k)$ is a maximum and therefore $p_\alpha$ peaks on $T_1$ at $z_{k_0}$, that is, $\|p_\alpha\|_{T_1} = |p_\alpha(z_{k_0})| > |p_\alpha(z)|$ for all $z \in T_1 \setminus \{z_{k_0}\}$.

If $\|p_\alpha\|_{T_2} < \|p_\alpha\|_{T_1}$, then we are done since then $p_\alpha$ peaks on $S = T_1 \cup T_2$ at $z_{k_0}$. If $\|p_\alpha\|_{T_2} \geq \|p_\alpha\|_{T_1}$, then let $M = \max\{|z_j| : j \in T_2\}$ and note that $0 \leq M < R$ by the definition of $T_2$. Thus, since $R > 1$ and $\|p_\alpha\|_{T_1} > 0$, for all $\alpha \in (0,1)$ with $\frac{\alpha-1}{\alpha} \neq z_{k_0}$, there exists a $K \in \mathbb{N}$ such that for all $m \geq K$,

$$M^m \|p_\alpha\|_{T_2} < R^m \|p_\alpha\|_{T_1}.$$

Clearly $p_{m,\alpha}(z) = z^m p_\alpha(z) = z^m(\alpha z + 1 - \alpha)$ is a convex polynomial and since $|p_{m,\alpha}(z)| = R^m |p_\alpha(z)|$ for $z \in T_1$ and $p_\alpha$ peaks on $T_1$ at $z_{k_0}$,

$$R^m \|p_\alpha\|_{T_1} = R^m |p_\alpha(z_{k_0})| = |p_{m,\alpha}(z_{k_0})|.$$

For $z \in T_2$, we have

$$|p_{m,\alpha}(z)| \leq M^m \|p_\alpha\|_{T_2} < R^m \|p_\alpha\|_{T_1} = |p_{m,\alpha}(z_{k_0})|.$$

Thus $p_{m,\alpha}$ peaks on $S = T_1 \cup T_2$ at $z_{k_0}$. It follows that

$$\|p_{m,\alpha}\|_S = |p_{m,\alpha}(z_{k_0})| = |z_{k_0}|^m |p_\alpha(z_{k_0})| = R^m \max_{|z_k|=R} |\alpha z_k + (1-\alpha)|.$$

Finally, notice that if $z_k \in S \setminus \mathbb{R}$, then $p_{m,\alpha}(z_k) = \alpha z_k^{m+1} + (1-\alpha) z_k^m$, thus $p_{m,\alpha}(z_k)$ is a convex-combination of $z_k^m$ and $z_k^{m+1}$ and since $z_k \notin \mathbb{R}$, at most one of $z_k^m$ and $z_k^{m+1}$ is a real number. It follows that at most one point on the line segment between $z_k^m$ and $z_k^{m+1}$ can be real. Thus there is at most one $\alpha_k \in (0,1)$ such that $p_{m,\alpha_k}(z_k)$ is a real number. Thus if $\alpha \in (0,1) \setminus \{\alpha_1, \ldots, \alpha_N\}$, then $p_{m,\alpha}(z_k) \notin \mathbb{R}$ for all $k$. □

**Lemma 4.2.** ***A One Variable Growth Lemma.*** *Let $\{M_n\}$ be a sequence of real numbers with $\lim M_n = +\infty$, $\{\varepsilon_n\}$ a sequence of complex numbers with $\lim |\varepsilon_n| = 0$, and $w$ a nonzero complex number. If $\theta$ is not an integer multiple of $\pi$, then there exists a subsequence $\{n_k\}$ of positive integers such that*

$$\lim_{k\to\infty} M_{n_k} Re\left(\varepsilon_{n_k} + e^{in_k\theta} w\right) = +\infty.$$

*Proof.* Writing $w$ as $|w|e^{i\alpha}$, we see that $e^{in\theta} w = |w| e^{i(n\theta+\alpha)}$. Since $\theta$ is not an integer multiple of $\pi$, either $\theta$ is an irrational multiple of $\pi$ or $\theta$ is a non-integer rational multiple of $\pi$. In the first instance, Kronecker's theorem implies that $\{e^{i(n\theta+\alpha)}\}$ is dense in the unit circle. Thus there is a subsequence $\{n_k\}$ such that $Re(e^{i(n_k\theta+\alpha)}) > \frac{1}{2}$. If $\theta$ is a non-integer rational multiple $\frac{p}{q}$ of $\pi$, then the points $e^{i(n\frac{p}{q}\pi+\alpha)}$ for $n = 1, 2, \cdots, 2q$ are evenly distributed about the unit circle. Hence there is an integer $n_0 \leq 2q$ such that $0 < Re\left(e^{i(n_0\frac{p}{q}\pi+\alpha)}\right) \leq 1$. Letting $n_k = n_0 + 2qk$ we see that, for all $k$,

$$0 < Re(e^{i(n_k\frac{p}{q}\pi+\alpha)}) = Re(e^{i(n_0\frac{p}{q}\pi+\alpha)}) < 1.$$

Hence in all cases, we can assert that there exists a subsequence $\{n_k\}$ and a fixed $\delta > 0$ such that $Re(|w|e^{i(n_k\theta+\alpha)}) > \delta$. Since $\lim |\varepsilon_n| = 0$ and $\lim M_n = +\infty$, for a given $R > 0$, there exists an $N \in \mathbb{N}$ such that for all $k > N$ we have $|\varepsilon_{n_k}| < \delta/2$ and $M_{n_k} > R$. Thus,

$$M_{n_k} Re\left(\varepsilon_{n_k} + e^{in_k\theta} w\right) > R(-\delta/2 + \delta) > R\delta/2.$$

The lemma is now immediate as we may choose $R$ to be arbitrarily large. □

## 5. Convex-Cyclic Diagonal Matrices

**Theorem 5.1.** ***Diagonal Matrices.***

***Complex Case:***

*If $T = diag(\lambda_1, \lambda_2, \ldots, \lambda_N)$ is a diagonal matrix on $\mathbb{C}^N$, then $T$ is convex-cyclic if and only if the following hold:*

1. *the diagonal entries $\{\lambda_k\}_{k=1}^N$ are distinct;*
2. *$|\lambda_k| > 1$ for all $1 \leq k \leq N$;*
3. *$\lambda_j \neq \overline{\lambda_k}$ for all $1 \leq j, k \leq N$.*

***Real Case:***

*If $T = diag(\lambda_1, \lambda_2, \ldots, \lambda_N)$ is a diagonal matrix on $\mathbb{R}^N$, then $T$ is convex-cyclic if and only if the following hold:*

1. *the diagonal entries $\{\lambda_k\}_{k=1}^N$ are distinct;*
2. *$\lambda_k < -1$ for all $1 \leq k \leq N$;*

*Furthermore, in both cases, the convex-cyclic vectors for $T$ are precisely those vectors $\vec{v}$ for which every coordinate of $\vec{v}$ is non-zero; and such vectors are dense in $\mathbb{C}^N$ or $\mathbb{R}^N$.*

Notice that condition (3) above says that the eigenvalues of $T$ cannot come in conjugate pairs and that none of them can be real numbers.

*Proof.* **The Complex Case:** Let $\vec{v} = (1, 1, \ldots, 1) \in \mathbb{C}^N$. We begin by showing that $\vec{v}$ is a convex-cyclic vector for $T$ under the stated assumptions. According to the Hahn-Banach Criterion we must show that

$$\sup_{p \in \mathcal{CP}} Re \left\langle p(T)\vec{v}, \vec{f} \right\rangle = \infty$$

for every nonzero $\vec{f} = (f_1, f_2, \ldots, f_N) \in \mathbb{C}^N$. Notice that

$$Re\langle p(T)\vec{v}, \vec{f}\rangle = Re\Big(p(\lambda_1)\overline{f_1} + p(\lambda_2)\overline{f_2} + \cdots + p(\lambda_N)\overline{f_N}\Big) = \tag{1}$$

$$= Re \sum_{k=1}^{N} p(\lambda_k)\overline{f_k} = Re \sum_{k \in A} p(\lambda_k)\overline{f_k}$$

where $A = \{k : f_k \neq 0\}$. Now by our hypothesis we see that the subset $\{\lambda_k : k \in A\}$ of the eigenvalues satisfies the hypothesis of Theorem 4.1, thus there is a convex-polynomial $p$ such that $p$ peaks on the set $\{\lambda_k : k \in A\}$ at a point $\lambda_j$ where $j \in A$, and $p$ satisfies $m := |p(\lambda_j)| > 1$, and $p(\lambda_j)$ is not a real number.

Now consider the sequence of convex-polynomials $\{p(z)^n\}_{n=1}^\infty$. Referring to (1) and writing $\frac{p(\lambda_j)}{|p(\lambda_j|} = e^{i\theta}$ where $\theta$ is not a multiple of $\pi$, we have

$$Re\langle p(T)^n \vec{v}, \vec{f}\rangle = Re \sum_{k \in A} p(\lambda_k)^n \cdot \overline{f_k} = m^n \cdot Re\left[\sum_{k \in A}\left(\frac{p(\lambda_k)}{m}\right)^n \cdot \overline{f_k}\right] =$$

$$= m^n \cdot Re\left[e^{in\theta}\overline{f_j} + \sum_{k \in A, k \neq j}\left(\frac{p(\lambda_k)}{m}\right)^n \cdot \overline{f_k}\right] = m^n Re[e^{in\theta}\overline{f_j} + \varepsilon_n]$$

where $\varepsilon_n \to 0$ since $\left|\frac{p(\lambda_k)}{m}\right| < 1$ for all $k \in A$, $k \neq j$. Since $m > 1$, $m^n \to \infty$. Since $\theta$ is not a multiple of $\pi$, Lemma 4.2 with $M_n = m_n$, $\varepsilon_n = \sum_{k\in A, k\neq j} \left(\frac{p(\lambda_k)}{m}\right)^n \cdot \overline{f_k}$, and $w = \overline{f_j}$, implies that

$$\sup_{n\geq 1} m^n Re[e^{in\theta}\overline{f_j} + \varepsilon_n] = \infty$$

and thus $\sup_{n\geq 1} Re\langle p(T)^n \vec{v}, \vec{f}\rangle = \infty$ as desired. It now follows that $T$ is convex-cyclic with convex-cyclic vector $\vec{v} = (1, 1, \ldots, 1)$.

**The Real Case:**

The proof is essentially the same as the complex case but with the simplification that Theorem 4.1 and Lemma 4.2 are not needed. Clearly $\mathbb{C}$ is everywhere replaced with $\mathbb{R}$. With $A = \{k : f_k \neq 0\}$, the subset $\{|\lambda_k| : k \in A\}$ of the eigenvalues has a unique maximum at some $\lambda_j$. Thus the convex polynomial $p(x) = x$ peaks at $\lambda_j$ and $m := |p(\lambda_j)| = |\lambda_j| > 1$. Hence $\langle T^n\vec{v}, \vec{f}\rangle = \sum_{k\in A} \lambda_k^n \cdot f_k =$

$$= m^n \cdot \sum_{k\in A} \left(\frac{\lambda_k}{m}\right)^n \cdot f_k = m^n(-1)^n f_j + \sum_{k\in A, k\neq j} \left(\frac{\lambda_k}{m}\right)^n \cdot f_k.$$

Choosing $n_k$, all even or all odd, such that $(-1)^{n_k} f_j > 0$ for all $k$ and noting that each of the terms $\left(\frac{\lambda_k}{m}\right)^n \cdot f_k$ goes to zero as $n \to \infty$, we see that $\sup_{n\geq 1}\langle T^n\vec{v}, \vec{f}\rangle = \infty$ which implies $T$ is convex-cyclic. The remainder of the proof is identical.

*The Convex-Cyclic Vectors:* To describe the convex-cyclic vectors for $T$, in both the real and complex cases, it is clear that every component of a convex-cyclic vector must be non-zero. For the converse, let $D$ be any diagonal invertible matrix. Then $D$ commutes with $T$ and has dense range (in fact it's onto), from this it follows that since $\vec{v} = (1, 1, \ldots, 1)$ is a convex-cyclic vector for $T$, then $Dv$ is also a convex-cyclic vector for $T$. Since $D$ can be any invertible diagonal matrix, it follows that $Dv$ can be any vector all of whose coordinates are non-zero. Thus all such vectors are convex-cyclic vectors for $T$. □

The following corollary is the case where all the diagonal entries in the diagonal matrix $T$ in Theorem 5.1 have absolute value equal to $r$, but they cannot be real or complex conjugates of one another. This is the difficult case in proving Theorem 5.1 which we were able to avoid by making use of peaking convex-polynomials. The $N = 2$ case of the next lemma is the One-Variable Growth Lemma (see Lemma 4.2) which uses Kronecker's Theorem.

**Corollary 5.2.** ***A Multivariable Growth Lemma.*** *If $\{f_k\}_{k=1}^N$ are complex numbers, not all zero, $r > 1$, $\{\theta_k\}_{k=1}^N$ are real numbers satisfying $\theta_i \neq \pm\theta_j$ (mod $2\pi$) and $\theta_j \neq n\pi$ for $n \in \mathbb{Z}$ and all $1 \leq j \leq N$, then*

$$\sup_{n\geq 1} r^n \cdot Re\left(\sum_{k=1}^N e^{in\theta_k} f_k\right) = \infty.$$

*Proof.* Let $T$ be the diagonal matrix with $\lambda_k = re^{i\theta_k}$ as its $k^{th}$ diagonal entry. Our hypothesis tells us that the $\{\lambda_k\}_{k=1}^N$ are distinct, have absolute value greater than one, and no two of them are complex conjugates of each other and none of them are real. Thus Theorem 5.1 applies to say that $T$ is convex-cyclic and as such with $\vec{f} = (f_1, \ldots, f_N) \neq \vec{0}$ we must have $\sup_{n\geq 1}\langle T^n\vec{v}, \vec{f}\rangle = \infty$ where $\vec{v} = (1, 1, \ldots, 1)$. □

## 6. Interpolating Convex-Polynomials

**Lemma 6.1.** ***Dense Convex Sets***

*If $C$ is a dense convex set in a finite dimensional real or complex vector space $X$, then $C = X$.*

*Proof.* Here is the idea of the proof in $\mathbb{R}^n$. Let $\vec{v} = (v_1, v_2, \ldots, v_n) \in \mathbb{R}^n$. Clearly, $\vec{v}$ lies in the interior of a sufficiently large $n$-cell (a product of intervals) in $\mathbb{R}^n$. Now since $C$ is dense in $\mathbb{R}^n$ we may approximate each of the $2^n$ vertices of this $n$-cell sufficiently close by vectors $\{\vec{x}_k : 1 \le k \le 2^n\}$ in $C$ so that $\vec{v}$ lies in the interior of the convex hull of $\{\vec{x}_k : 1 \le k \le 2^n\}$. However since $C$ is convex we have that $\vec{v} \in co(\{\vec{x}_k\}) \subseteq C$. Since $\vec{v}$ was arbitrary in $\mathbb{R}^n$ we have that $C = \mathbb{R}^n$. A similar proof works in $\mathbb{C}^n$ and every real or complex finite dimensional vector space is isomorphic to $\mathbb{R}^n$ or $\mathbb{C}^n$, so the Lemma follows. □

**Theorem 6.2.** ***Interpolating Prescribed Values.***
***Complex Case:*** *If $\{z_k\}_{k=1}^n$ are distinct complex numbers satisfying*

(1) $|z_k| > 1$ *for all* $1 \le k \le n$;
(2) $Im(z_k) \neq 0$ *for all* $1 \le k \le n$; *and*
(3) $z_k \neq \overline{z}_j$ *for all* $1 \le j, k \le n$,

*then for any finite set of complex numbers $\{w_k\}_{k=1}^n$, there exists a convex-polynomial $p$ such that $p(z_k) = w_k$ for all $1 \le k \le n$.*

***Real Case:*** *If $\{x_k\}_{k=1}^n$ are distinct real numbers satisfying $x_k < -1$ for all $1 \le k \le n$ and $\{y_k\}_{k=1}^n$ is a finite set of real numbers, then there exists a convex-polynomial $p$ such that $p(x_k) = y_k$ for all $1 \le k \le n$.*

*Proof. The Complex Case.* Let $C$ be the set of all vectors $(w_1, w_2, \ldots, w_n) \in \mathbb{C}^n$ such that there exists a convex-polynomial $p$ satisfying $p(z_k) = w_k$ for all $1 \le k \le n$. Then $C$ is a convex subset of $\mathbb{C}^n$. We want to show that $C = \mathbb{C}^n$. Let $T$ be the diagonal matrix with diagonal entries $(z_1, z_2, \ldots, z_n)$. Then conditions $(1), (2)$, and $(3)$ in our hypothesis together with Theorem 5.1 imply that the matrix $T$ is convex-cyclic and the vector $\vec{v} = (1, 1, \ldots, 1)$ is a convex cyclic vector for $T$. It follows that $\{p(T)\vec{v} : p \in \mathcal{CP}\}$ is dense in $\mathbb{C}^n$. Since $p(T)\vec{v} = \big(p(z_1), p(z_2), \ldots, p(z_n)\big)$ we see that the set $C$ of values which can be interpolated is a dense subset of $\mathbb{C}^n$. However, by Lemma 6.1 the dense convex subset $C$ of $\mathbb{C}^n$ is equal to $\mathbb{C}^n$, as desired. The proof of the real case is similar. □

## 7. Convex-Cyclic Jordan Matrices

Recall that the $k \times k$ lower Jordan block with eigenvalue $\lambda$, denoted by $J_k(\lambda)$, is a $k \times k$ matrix with $\lambda$ along the main diagonal and ones along the subdiagonal, and zeros elsewhere. Below is $J_4(\lambda)$.

$$J_4(\lambda) = \begin{bmatrix} \lambda & 0 & 0 & 0 \\ 1 & \lambda & 0 & 0 \\ 0 & 1 & \lambda & 0 \\ 0 & 0 & 1 & \lambda \end{bmatrix}.$$

Powers of the matrix $J_k(\lambda)$ follow a simple pattern:

$$J_4(\lambda)^n = \begin{bmatrix} \lambda^n & 0 & 0 & 0 \\ n\lambda^{n-1} & \lambda^n & 0 & 0 \\ \frac{n(n-1)}{2}\lambda^{n-2} & n\lambda^{n-1} & \lambda^n & 0 \\ \frac{n(n-1)(n-2)}{3!}\lambda^{n-3} & \frac{n(n-1)}{2}\lambda^{n-2} & n\lambda^{n-1} & \lambda^n \end{bmatrix} =$$

$$\begin{bmatrix} \lambda^n & 0 & 0 & 0 \\ \binom{n}{1}\lambda^{n-1} & \lambda^n & 0 & 0 \\ \binom{n}{2}\lambda^{n-2} & \binom{n}{1}\lambda^{n-1} & \lambda^n & 0 \\ \binom{n}{3}\lambda^{n-3} & \binom{n}{2}\lambda^{n-2} & \binom{n}{1}\lambda^{n-1} & \lambda^n \end{bmatrix}$$

The terms down the first column are simply $\frac{1}{k!}p^{(k)}(\lambda)$ where $p(z) = z^n$ and the matrix is a lower-triangular Toeplitz matrix (constant along the diagonals). Note also that the coefficient of $\lambda^{n-j}$ is $\binom{n}{j}$. In fact if $p$ is any polynomial (or analytic function), then the same pattern applies.

**Proposition 7.1.** *If $n \geq 1$ and $p$ is a polynomial, then $p(J_n(\lambda))$ is a lower triangular Toeplitz matrix of the following form:*

$$p(J_n(\lambda)) = \begin{bmatrix} p(\lambda) & 0 & \dots & 0 & 0 & 0 \\ p'(\lambda) & p(\lambda) & 0 & \dots & 0 & 0 \\ \frac{p^{(2)}(\lambda)}{2!} & p'(\lambda) & p(\lambda) & 0 & \dots & 0 \\ \vdots & \vdots & p'(\lambda) & \ddots & 0 & 0 \\ \frac{p^{(n-2)}(\lambda)}{(n-2)!} & \frac{p^{(n-3)}(\lambda)}{(n-3)!} & \dots & p'(\lambda) & p(\lambda) & 0 \\ \frac{p^{(n-1)}(\lambda)}{(n-1)!} & \frac{p^{(n-2)}(\lambda)}{(n-2)!} & \dots & \frac{p^{(2)}(\lambda)}{2!} & p'(\lambda) & p(\lambda) \end{bmatrix}. \tag{1}$$

*In particular, $p(J_n(\lambda)) = 0$ if and only if $p^{(j)}(\lambda) = 0$ for all $0 \leq j \leq (n-1)$; in other words, if and only if $p$ has a zero of order $n$ at $\lambda$.*

Every complex $n \times n$ matrix $T$ is similar to a Jordan matrix which is a direct sum of Jordan blocks $J_k(\lambda)$ of various sizes. Given a positive integer $N$, let $e_{N,k}$ be the unit basis vector of length $N$ with a one in the $k^{th}$ position.

**Proposition 7.2. Convex-Cyclic Jordan blocks.**

***The Complex Case:*** *For $\lambda \in \mathbb{C}$ and $m \geq 2$ the Jordan block $J_m(\lambda)$ is convex-cyclic on $\mathbb{C}^m$ if and only if $|\lambda| > 1$ and $\lambda \notin \mathbb{R}$. Furthermore, a vector $\vec{v} = (v_1, v_2, \dots, v_m) \in \mathbb{C}^m$ is a convex-cyclic vector for $J_m(\lambda)$ if and only if $\vec{v}$ is a cyclic vector for $J_m(\lambda)$, which holds, if and only if $v_1 \neq 0$.*

***The Real Case:*** *For $\lambda \in \mathbb{R}$ and $m \geq 2$ the Jordan block $J_m(\lambda)$ is convex-cyclic on $\mathbb{R}^m$ if and only if $\lambda < -1$. Furthermore, a vector $\vec{v} = (v_1, v_2, \dots, v_m) \in \mathbb{R}^m$ is a convex-cyclic vector for $J_m(\lambda)$ if and only if $\vec{v}$ is a cyclic vector for $J_m(\lambda)$, which holds if and only if $v_1 \neq 0$.*

*Proof.* **The Complex Case.** If $J_m(\lambda)$ is convex-cyclic, then by Proposition 2.3 we know that $|\lambda| > 1$ and $\lambda \notin \mathbb{R}$. Now suppose that $|\lambda| > 1$ and $\lambda \notin \mathbb{R}$ and we will show that $J_m(\lambda)$ is convex-cyclic. Let $\vec{e_1} = (1, 0, 0, \dots, 0) \in \mathbb{C}^m$. We'll begin by showing that $\vec{e_1}$ is a convex-cyclic vector for $J_m(\lambda)$ and then afterwards show that other vectors are also convex-cyclic vectors. First write $\lambda = re^{i\theta}$ where $r > 1$ and $\theta$ is not a multiple of $\pi$. Let $\vec{f} = (f_1, f_2, \dots, f_m) \in \mathbb{C}^m \setminus \{\vec{0}\}$. We must show that

$$\sup_{n \geq 1} Re \left\langle J_m(\lambda)^n \vec{e_1}, \vec{f} \right\rangle = \infty.$$

Let $j$ be the largest integer such that $f_j \neq 0$ and we will factor out $\binom{n}{j-1} r^{n-(j-1)}$ in the expression below.

Note that for $n \geq 2m$ the sequence $\left\{\binom{n}{k}\right\}_{k=0}^{m}$ is increasing (the rows of Pascal's triangle first increase, then decrease), thus for large $n$ ($n \geq 2j$) we have

$$Re\left\langle J_m(\lambda)^n \vec{e_1}, \vec{f}\right\rangle = Re\sum_{k=0}^{j-1}\binom{n}{k}\lambda^{n-k}\cdot\overline{f}_{k+1} =$$

$$Re\left[\binom{n}{0}\lambda^n\cdot\overline{f}_1 + \binom{n}{1}\lambda^{n-1}\cdot\overline{f}_2 + \cdots + \binom{n}{j-2}\lambda^{n-(j-2)}\cdot\overline{f}_{j-1} + \binom{n}{j-1}\lambda^{n-(j-1)}\cdot\overline{f}_j\right]$$

$$= \binom{n}{j-1}r^{n-(j-1)}\cdot Re\left[\left(\sum_{k=0}^{j-2}\frac{\binom{n}{k}}{\binom{n}{j-1}}r^{j-k-1}e^{i(n-k)\theta}\cdot\overline{f}_{k+1}\right) + e^{i(n-(j-1))\theta}\cdot\overline{f}_j\right]$$

Note that $\frac{\binom{n}{k}}{\binom{n}{j-1}} \to 0$ as $n \to \infty$ in the sum above since $0 \leq k < j-1 < n/2$. Also, $\binom{n}{j-1}r^{n-(j-1)} \to \infty$ as $n \to \infty$ since $r > 1$. Letting $M_n = \binom{n}{j-1}r^{n-(j-1)}$, $\varepsilon_n = \sum_{k=0}^{j-2}\frac{\binom{n}{k}}{\binom{n}{j-1}}r^{j-k-1}e^{i(n-k)\theta}\cdot\overline{f}_{k+1}$, and $w = \overline{f}_j$, and noting that $\theta$ is not a multiple of $\pi$, Lemma 4.2 implies that there is a subsequence $\{n_k\}$ such that

$$\lim_{k\to\infty} Re\left\langle J_m(\lambda)^{n_k}\vec{e_1}, \vec{f}\right\rangle = \infty.$$

It follows that $J_m(\lambda)$ is convex-cyclic with convex-cyclic vector $\vec{e_1} = (1, 0, 0, \ldots, 0)$.

In order to prove the last claim of the theorem, let $\vec{v} = (v_1, v_2, \ldots, v_m) \in \mathbb{C}^m$ with $v_1 \neq 0$ and we will show that $\vec{v}$ is a convex-cyclic vector for $J_m(\lambda)$. Let $p$ be the polynomial $p(z) = \sum_{k=0}^{m-1} v_{k+1}z^k$ and consider the finite Toeplitz matrices $T$ and $S$ given below:

$$T = \begin{bmatrix} v_1 & 0 & 0 & \cdots & 0 \\ v_2 & v_1 & 0 & \ddots & \vdots \\ v_3 & v_2 & \ddots & \ddots & 0 \\ \vdots & \ddots & \ddots & \ddots & 0 \\ v_m & \ddots & v_3 & v_2 & v_1 \end{bmatrix} \qquad S = \begin{bmatrix} 0 & 0 & 0 & \cdots & 0 \\ 1 & 0 & 0 & \ddots & \vdots \\ 0 & 1 & \ddots & \ddots & 0 \\ \vdots & \ddots & \ddots & \ddots & 0 \\ 0 & \ddots & 0 & 1 & 0 \end{bmatrix}.$$

Note the finite Toeplitz matrices $T$ and $J_m(\lambda)$ are both polynomials in $S$, so they commute. Also, $T$ is invertible since $v_1 \neq 0$.

We showed above that $\vec{e_1} = (1, 0, 0, \ldots, 0) \in \mathbb{C}^m$ is a convex-cyclic vector for $J_m(\lambda)$. Hence $T\vec{e_1} = (v_1, v_2, \ldots, v_m)$ is also a convex-cyclic vector for $J_m(\lambda)$ since $p(J_m(\lambda))T\vec{e_1} = Tp(J_m(\lambda))\vec{e_1}$ for any convex-polynomial $p$ and $T$ being invertible means $T$ maps a dense set to a dense set.

Conversely, if $\vec{v} = (0, v_2, \ldots, v_m) \in \mathbb{C}^m$, then from (1) we see that for any polynomial $p$, $p(J_m(\lambda))\vec{v} = (0, w_2, w_3, \ldots, w_m)$ for some scalars $\{w_k\}_{k=2}^m$. Thus $\{p(J_m(\lambda))\vec{v} : p \text{ is a polynomial}\}$ cannot be dense in $\mathbb{C}^m$, thus $\vec{v}$ is not a cyclic vector for $J_m(\lambda)$.

**The Real Case.** This is entirely similar to the complex case with $\mathbb{C}$ replaced by $\mathbb{R}$ throughout and $\lambda = -r$ where $r > 1$. Lemma 3.2 is not required as

$$\left\langle J_m(\lambda)^n\vec{e_1}, \vec{f}\right\rangle = \sum_{k=0}^{j-1}\binom{n}{k}(-r)^{n-k}\cdot f_{k+1} =$$

$$= \binom{n}{j-1}r^{n-(j-1)}\cdot\left[\left(\sum_{k=0}^{j-2}\frac{\binom{n}{k}}{\binom{n}{j-1}}r^{j-k-1}(-1)^{(n-k)}\cdot f_{k+1}\right) + (-1)^{(n-(j-1))}\cdot f_j\right].$$

As above $\frac{\binom{n}{k}}{\binom{n}{j-1}} \to 0$ as $n \to \infty$ since $0 \le k < j-1 < n/2$. So $\left\langle J_m(\lambda)^n \vec{e_1}, \vec{f} \right\rangle \to \infty$ for a subsequence of $\{n\}_{n=1}^{\infty}$ chosen so that $(-1)^{(n-(j-1))} \cdot f_j = |f_j| > 0$. The remainder of the proof is exactly the same as the complex case. □

The next result examines a direct sum of a diagonal matrix and an $m \times m$ Jordan block.

**Theorem 7.3.** ***Diagonal Matrix direct sum a Jordan Block.***

***Complex Case:*** *Suppose $T = D \oplus J$ where $D = diag(\lambda_1, \lambda_2, \dots, \lambda_N)$ is a diagonal matrix on $\mathbb{C}^N$ and $J = J_m(\lambda_{N+1})$ is an $m \times m$ Jordan block, $m \ge 1$, with eigenvalue $\lambda_{N+1}$, then $T$ is convex-cyclic if and only if the following hold:*

1. *the eigenvalues $\{\lambda_k\}_{k=1}^{N+1}$ are distinct;*
2. *$|\lambda_k| > 1$ for all $1 \le k \le N+1$;*
3. *$\lambda_j \ne \overline{\lambda_k}$ for all $1 \le j, k \le N+1$.*

***Real Case:*** *Suppose $T = D \oplus J$ where $D = diag(\lambda_1, \lambda_2, \dots, \lambda_N)$ is a diagonal matrix on $\mathbb{R}^N$ and $J = J_m(\lambda_{N+1})$ is an $m \times m$ Jordan block, $m \ge 1$, with eigenvalue $\lambda_{N+1}$, then $T$ is convex-cyclic if and only if the following hold:*

1. *the eigenvalues $\{\lambda_k\}_{k=1}^{N+1}$ are distinct;*
2. *$\lambda_k < -1$ for all $1 \le k \le N+1$.*

*Furthermore, in both the real and complex cases the convex-cyclic vectors for $T$ are precisely those vectors $\vec{w} = (u_1, u_2, \dots, u_N, v_1, v_2, \dots, v_m) = \vec{u} \oplus \vec{v}$ for which $u_k \ne 0$ for all $1 \le k \le N$ and $v_1 \ne 0$. That is, $\vec{w}$ is convex-cyclic for $T$ if and only if $\vec{u}$ is a cyclic vector for $D$ and $\vec{v}$ is a cyclic vector for $J$.*

*Proof.* **The Complex Case:** First note that $J$ is convex-cyclic by Proposition 7.2 and $D$ is convex cyclic by Theorem 5.1. We will apply Proposition 3.3, for which we need to find a convex-polynomial $p_0$ such that $p_0(D)$ is convex-cyclic and $p_0(J) = 0$. By Theorem 6.2, there is a convex-polynomial $q$ such that $q(\lambda_k) = 2ke^{i\sqrt{2}\cdot\pi}$ for $1 \le k \le N$ and $q(\lambda_{N+1}) = 0$. Let $p_0(z) = q(z)^m$. Then $|p_0(\lambda_k)| = 2^m k^m > 1$ for all $1 \le k \le N$, $p_0(\lambda_k) = 2^m k^m e^{im\sqrt{2}\cdot\pi}$ is not a real number, and clearly $p_0(\lambda_j)$ is not equal to $\overline{p_0(\lambda_k)}$ for $1 \le j, k \le N$ since they have different absolute values. It then follows from Theorem 5.1 that $p_0(D)$ is convex-cyclic. Since $q(z)$ has a zero at $\lambda_{N+1}$, $p_0(z) = q(z)^m$ has a zero of order $m$ at $\lambda_{N+1}$. Thus $p_0^{(j)}(\lambda_{N+1}) = 0$ for all $0 \le j \le (m-1)$. Since $J = J_m(\lambda_{N+1})$ is an $m \times m$ Jordan block, Proposition 7.1 implies $p_0(J) = 0$. It now follows from Proposition 3.3 that $T = D \oplus J$ is convex-cyclic.

It also follows from Proposition 3.3 that if $\vec{u}$ is a convex-cyclic vector for $p_0(D)$ and $\vec{v}$ is a convex-cyclic vector for $J$, then $\vec{u} \oplus \vec{v}$ is a convex-cyclic vector for $T$. From Theorem 5.1 we see that any vector $\vec{u}$ all of whose coordinates are non-zero is a convex-cyclic vector for $p_0(D)$. Also by Proposition 7.2 we see that any vector $\vec{v}$ whose first coordinate is non-zero is a convex-cyclic vector for $J$. Thus such vectors $\vec{u} \oplus \vec{v}$ are convex-cyclic for $T$. These conditions are also clearly necessary for a vector to be convex-cyclic for $T$.

**The Real Case:** The proof is similar to the complex case, except now choose $q$ so that $q(\lambda_k) = -2k$ for $1 \le k \le N$ and $q(\lambda_{N+1}) = 0$. Then choose an odd integer $r \ge m$ and let $p_0(x) = q(x)^r$. Then $\{p_0(\lambda_k)\}_{k=1}^{N}$ are distinct and $p_0(\lambda_k) = (-2k)^r < -1$ for all $1 \le k \le N$, which implies that $p_0(D)$ is convex-cyclic (see Theorem 5.1). Also $p_0$ has a zero of order at least $r$ $(\ge m)$ at $\lambda_{N+1}$ so $p_0(J) = 0$ by Proposition 7.1 . It now follows from Proposition 3.3 that $T$ is convex-cyclic. The rest of the proof is the same as in the complex case. □

**Theorem 7.4.** ***Interpolating values & Derivatives at one point.***

***Complex Case:*** *If $\{z_k\}_{k=1}^{n+1}$ are distinct complex numbers satisfying*

(1) $|z_k| > 1$ *for all* $1 \le k \le (n+1)$;
(2) $Im(z_k) \neq 0$ *for all* $1 \le k \le (n+1)$; *and*
(3) $z_j \neq \overline{z}_k$ *for all* $1 \le j, k \le (n+1)$,

*then for any finite set of complex numbers* $\{w_{0,k}\}_{k=1}^{n+1} \cup \{w_{j,n+1}\}_{j=1}^{m}$ *there exists a convex-polynomial* $p$ *such that* $p(z_k) = p^{(0)}(z_k) = w_{0,k}$ *for all* $1 \le k \le n+1$ *and* $p^{(j)}(z_{n+1}) = w_{j,n+1}$ *for* $1 \le j \le m$.

***Real Case:*** *If* $\{x_k\}_{k=1}^{n}$ *are distinct real numbers satisfying* $x_k < -1$ *for all* $1 \le k \le n$ *and* $\{y_k\}_{k=1}^{n} \cup \{y_{j,n}\}_{j=1}^{m}$ *is a finite set of real numbers, then there exists a convex-polynomial* $p$ *such that* $p(x_k) = y_k$ *for all* $1 \le k \le n$ *and* $p^{(j)}(x_n) = y_{j,n}$ *for* $1 \le j \le m$.

*Proof.* Let $C$ be the set of all vectors in $\mathbb{C}^{n+1+m}$ of the form

$$(w_{0,1}, w_{0,2}, \dots, w_{0,n+1}, w_{1,n+1}, w_{2,n+1}, \dots, w_{m,n+1})$$

such that there exists a convex-polynomial $p$ satisfying $p(z_k) = w_{0,k}$ for $1 \le k \le n+1$ and such that $p^{(j)}(z_{n+1}) = w_{j,n+1}$ for $1 \le j \le m$. Then $C$ is a convex subset of $\mathbb{C}^{n+1+m}$. We want to show that $C = \mathbb{C}^{n+1+m}$. Let $D$ be the $(n+1) \times (n+1)$ diagonal matrix with diagonal entries $(z_1, z_2, \dots, z_{n+1})$ and let $J_m(z_{n+1})$ be the $m \times m$ Jordan block with eigenvalue $z_{n+1}$. Since conditions $(1), (2)$, and $(3)$ hold, then by Theorem 7.3 the matrix $T = D \oplus J_m(z_{n+1})$ is convex-cyclic and the vector $\vec{v} = (\vec{v}_1, \vec{v}_2) = (1, 1, \dots, 1, 0, 0, \dots, 0)$ where $\vec{v}_1 = (1, \dots, 1) \in \mathbb{C}^{n+1}$ and $\vec{v}_2 = (1, 0, 0, \dots, 0) \in \mathbb{C}^m$ is a convex cyclic vector for $T$. It follows that $\{p(T)\vec{v} : p \in \mathcal{CP}\}$ is dense in $\mathbb{C}^{n+1+m}$. Since $p(T)\vec{v} = \big(p(z_1), p(z_2), \dots, p(z_{n+1}), p'(z_{n+1}), \dots, p^{(m)}(z_{n+1})\big)$ we see that the set $C$ of values which can be interpolated is a dense subset of $\mathbb{C}^{n+1+m}$.

However, by Lemma 6.1 the dense convex subset $C$ of $\mathbb{C}^{n+1+m}$ is equal to $\mathbb{C}^{n+1+m}$, as desired. The proof of the real case is similar. □

With the above interpolation result, we are now prepared to prove exactly which Jordan matrices are convex-cyclic.

**Theorem 7.5.** ***A Jordan Matrix.***

***The Complex Case:***

*If* $J = \bigoplus_{k=1}^{N} J_{n_k}(\lambda_k)$ *is a Jordan matrix on* $\mathbb{C}^p$ *where* $p = \sum_{k=1}^{N} n_k$, *then* $J$ *is convex-cyclic on* $\mathbb{C}^p$ *if and only if the following hold:*

(1) *the eigenvalues* $\{\lambda_k\}_{k=1}^{N}$ *are distinct;*
(2) $|\lambda_k| > 1$ *for all* $1 \le k \le N$;
(3) $\lambda_j \neq \overline{\lambda_k}$ *for all* $1 \le j, k \le N$.

***The Real Case:***

*If* $J = \bigoplus_{k=1}^{N} J_{n_k}(\lambda_k)$ *is a Jordan matrix on* $\mathbb{R}^p$ *where* $p = \sum_{k=1}^{N} n_k$, *then* $J$ *is convex-cyclic on* $\mathbb{R}^p$ *if and only if the eigenvalues* $\{\lambda_k\}_{k=1}^{N}$ *are distinct and* $\lambda_k < -1$ *for all* $1 \le k \le N$.

*Furthermore, in both the real and the complex cases, the convex-cyclic vectors for* $J$ *are precisely those vectors of the form* $\vec{v} = (\vec{v}_1, \vec{v}_2, \dots, \vec{v}_N)$ *where for all* $1 \le k \le N$, $\vec{v}_k$ *is a convex-cyclic vector for* $J_{n_k}(\lambda_k)$; *that is,* $\vec{v}_k \in \mathbb{F}^{n_k}$ *and the first coordinate* $\vec{v}_k(1)$ *of* $\vec{v}_k$ *must be non-zero; where* $\mathbb{F}$ *equals* $\mathbb{R}$ *or* $\mathbb{C}$.

*Proof.* **Complex Case.** We shall use Proposition 3.3 and induction to prove this theorem. We shall do induction on the number $N$ of Jordan blocks that appear in $J$. If $N = 1$, then Proposition 7.2 applies and says that $J$ is convex-cyclic and any vector whose first coordinate is non-zero is a convex-cyclic vector.

Now suppose that the theorem holds for any $N$ Jordan blocks; in other words, suppose that the direct sum of any $N$ Jordan blocks satisfying conditions (1), (2), and (3) is convex-cyclic and any vector $\vec{v}$ as described in the theorem is a convex-cyclic vector and we will show that the direct sum of any $N+1$ Jordan blocks satisfying (1), (2), and (3) is convex-cyclic and has the specified set of convex-cyclic vectors.

So, let $\{J_{n_k}(\lambda_{n_k}) : 1 \leq k \leq N+1\}$ be a collection of $(N+1)$ Jordan blocks and suppose that $\{\lambda_{n_k}\}_{k=1}^{N+1}$ satisfies the conditions (1), (2), and (3) of our hypothesis. Let $J = J_{n_1}(\lambda_1) \oplus J_{n_2}(\lambda_2) \oplus \cdots \oplus J_{n_{N+1}}(\lambda_{N+1})$. We will show that $J$ is convex-cyclic and that it has the stated convex-cyclic vectors.

Let $T_1 = J_{n_{N+1}}(\lambda_{N+1})$ and $T_2 = J_{n_1}(\lambda_1) \oplus J_{n_2}(\lambda_2) \oplus \cdots \oplus J_{n_N}(\lambda_N)$ and we'll apply Proposition 3.3. Since (1), (2) and (3) hold, by Theorem 7.4, there is a convex-polynomial $p_0$ such that $p_0(\lambda_k) = 0$ for all $1 \leq k \leq N$ and $p_0(\lambda_{N+1}) = 2i$, $p_0'(\lambda_{N+1}) = 1$ and $p_0^{(k)}(\lambda_{N+1}) = 0$ for $2 \leq k \leq n_{N+1}$. It then follows from Proposition 7.1 that $p_0(T_1) = J_{n_{N+1}}(2i)$ and $p_0(T_2)$ has zero as its only eigenvalue. Thus by Proposition 7.2 we know that $p_0(T_1)$ is convex-cyclic. Also since $p_0(T_2)$ has zero as its only eigenvalue it is nilpotent, hence certainly power bounded. Thus Proposition 3.3 applies to say that $J = T_2 \oplus T_1$ is convex-cyclic, as desired.

Also, if $\vec{v} = (\vec{v}_1, \vec{v}_2, \ldots, \vec{v}_{N+1})$ is a vector as stated in our hypothesis, then $\vec{w} = (\vec{v}_1, \vec{v}_2, \ldots, \vec{v}_N)$ is a convex-cyclic vector for $T_2$ by our induction hypothesis and $\vec{v}_{N+1}$ is a convex-cyclic vector for $T_1$ by Proposition 7.2, thus by Proposition 3.3 we have that $(\vec{w}, \vec{v}_{N+1})$ is a convex-cyclic vector for $J$. Thus $J$ has the stated set of convex-cyclic vectors. Conversely if $\vec{v} = (\vec{v}_1, \vec{v}_2, \ldots, \vec{v}_{N+1})$ is a convex-cyclic vector for $J$, then each $\vec{v}_k$ is a convex-cyclic vector for $J_{n_k}$ and by Proposition 7.2 we must have that the first coordinate of $\vec{v}_k$ must be nonzero. Thus $J$ has the stated set of convex-cyclic vectors.

Finally, simply apply Proposition 2.3 in order to see that conditions (1), (2) and (3) are necessary.

**Real Case.** This case is naturally similar to the complex case. Use the real case of Theorem 7.4 to choose a convex-polynomial $p$ that satisfies $p(\lambda_k) = 0$ for all $1 \leq k \leq N$, $p(\lambda_{N+1}) = -2$, $p'(\lambda_{N+1}) = 1$ and $p^{(k)}(\lambda_{N+1}) = 0$ for $2 \leq k \leq n_{N+1}$. Then the only eigenvalue for $p(T_2)$ is 0 and thus $p(T_2)$ is nilpotent. Also, $p(T_1) = J_{n_{N+1}}(-2)$ which is convex-cyclic by Proposition 7.2. Thus Proposition 3.3 implies that $J = T_2 \oplus T_1$ is convex-cyclic. The rest of the proof is the same. □

Since the previous result gives a larger class of matrices that are convex-cyclic we get a stronger interpolation theorem.

**Theorem 7.6.** ***Convex-Polynomial Interpolation***
***The Complex Case:*** *Let $\{z_k\}_{k=1}^n$ be a finite set of complex numbers. Then the following are equivalent:*

*(a) for any finite set $\{w_{j,k} : 0 \leq j \leq N, 1 \leq k \leq n\}$ of complex numbers there exists a convex-polynomial $p$ such that $p^{(j)}(z_k) = w_{j,k}$ for all $0 \leq j \leq N$ and $1 \leq k \leq n$.*

*(b) The complex numbers $\{z_k\}_{k=1}^n$ are distinct, $\{z_k\}_{k=1}^n \subseteq \mathbb{C} \setminus \overline{\mathbb{D}}$ and $z_j \neq \overline{z}_k$ for all $1 \leq j, k \leq n$.*

***The Real Case:*** *If $\{x_k\}_{k=1}^n$ is a finite set of real numbers, then the following are equivalent:*

*(i) for any finite set $\{y_{j,k} : 0 \leq j \leq N, 1 \leq k \leq n\}$ of real numbers, there exists a convex-polynomial $p$ such that $p^{(j)}(x_k) = y_{j,k}$ for all $0 \leq j \leq N$ and $1 \leq k \leq n$.*

*(ii) The real numbers $\{x_k\}_{k=1}^n$ are distinct and satisfy $\{x_k\}_{k=1}^n \subseteq (-\infty, -1)$.*

*Proof.* Let's begin with the complex case and show that $(b)$ implies $(a)$. Let $C$ be the set of all $(N+1) \times n$ matrices $(w_{j,k})_{0 \leq j \leq N, 1 \leq k \leq n}$ with complex entries such that there exists a convex-polynomial $p$ satisfying $p^{(j)}(z_k) = w_{j,k}$ for all $0 \leq j \leq N$ and $1 \leq k \leq n$. In other words $C$ consists of all matrices of the form $(p^{(j)}(z_k))_{0\leq j\leq N, 1\leq k\leq n}$ where $p$ is a convex-polynomial. Then $C$ is a convex subset of $M_{N+1,n}(\mathbb{C})$, the vector space of all complex matrices of size $(N+1)\times n$.

Let $J = \bigoplus_{k=1}^{n} J_{N+1}(z_k)$ be the direct sum of $n$ Jordan blocks each having (the same) size $(N+1)\times(N+1)$ and the $k^{th}$ block having eigenvalue $z_k$. Since $(b)$ holds, we know from Theorem 7.5 that $J$ is convex-cyclic and that $\vec{v} = (\vec{e}_1, \vec{e}_1, \ldots, \vec{e}_1) \in \mathbb{C}^{n(N+1)}$ is a convex-cyclic vector for $J$ where $\vec{e}_1 = (1, 0, 0, \ldots, 0) \in \mathbb{C}^{N+1}$. It follows that $\{p(J)\vec{v} : p \in \mathcal{CP}\}$ is dense in $\mathbb{C}^{n(N+1)} \cong M_{N+1,n}(\mathbb{C})$. Here the isomorphism $\cong$ is given by mapping the first few entries of the vector to the first column of the matrix and the next few entries of the vector to the second column and so forth. Since $p(J)\vec{v} \cong (p^{(j)}(z_k))_{0\leq j\leq N, 1\leq k\leq n}$, it follows that the set $C$ of matrices in $M_{N+1,n}(\mathbb{C})$ that can be interpolated by a convex-polynomial is dense in $M_{N+1,n}(\mathbb{C})$. Since $C$ is both convex and dense in $M_{N+1,n}(\mathbb{C})$ it follows from Lemma 6.1 that $C$ must be equal to $M_{N+1,n}(\mathbb{C})$. The fact that $(b)$ implies $(a)$ now follows. The fact that $(a)$ implies $(b)$ follows from part (3) of Proposition 2.1. The real case is similar. □

By using the fact that every matrix is similar to its Jordan Canonical Form we get the following result from Theorem 7.5.

**Theorem 7.7.** ***Convex-Cyclicity of Complex Matrices.***
*If $A$ is an $n \times n$ matrix, then $A$ is convex-cyclic on $\mathbb{C}^n$ if and only if $A$ is cyclic and its eigenvalues $\{\lambda_k\}_{k=1}^{n}$ satisfy $|\lambda_k| > 1$ for all $1 \leq k \leq n$ and $\lambda_j \neq \overline{\lambda}_k$ for all $1 \leq j, k \leq n$.*

*Proof.* If $A$ is convex-cyclic on $\mathbb{C}^n$, then certainly $A$ is cyclic and Proposition 2.3 implies that the eigenvalues must have absolute value strictly larger than one, they cannot be real, and none of them can come in conjugate pairs. For the converse, assume the stated conditions hold and consider the Jordan form $J$ for $A$. Since $A$ is cyclic, then each eigenvalue has a geometric multiplicity of one, so that each eigenvalue appears in exactly one of the Jordan blocks in $J$. Thus the conditions of Theorem 7.5 are satisfied, so it follows that $J$ is convex-cyclic and since $A$ is similar to $J$, we must have that $A$ is also convex-cyclic. □

## 8. Real Convex-Cyclic Matrices

Our previous results on real matrices have all been about matrices on $\mathbb{R}^n$ with real eigenvalues. In this section we consider the case of matrices on $\mathbb{R}^n$ with real and complex eigenvalues. This case is actually the most general case of all the cases considered in this paper. This case uses the previous results proven about real matrices and complex matrices via the complexification map.

8.1. **Brief Review of Jordan Canonical Forms.** If $T$ is a real matrix, its eigenvalues may be complex and in that case the *real Jordan form* for $T$ is useful. The real Jordan form uses the Jordan blocks $J_k(\lambda)$ when $\lambda$ is real and some additional real blocks with complex eigenvalues. Let

$$C_1(r,\theta) = \begin{bmatrix} r\cos(\theta) & -r\sin(\theta) \\ r\sin(\theta) & r\cos(\theta) \end{bmatrix} = r\begin{bmatrix} \cos(\theta) & -\sin(\theta) \\ \sin(\theta) & \cos(\theta) \end{bmatrix} = rR(\theta).$$

Then $C_1(r,\theta)$ has complex eigenvalues $a \pm ib = r\cos(\theta) \pm ir\sin(\theta) = re^{\pm i\theta}$ and $R(\theta)$ is the matrix that rotates by an angle of $\theta$.

The $2k \times 2k$ real Jordan block $C_k(r,\theta)$ is the block lower-triangular matrix with $k$ copies of $C_1(r,\theta)$ down the main diagonal and with $2\times 2$ identity matrices on the block-subdiagonal. Below is an example:

$$C_3(r,\theta) = \begin{bmatrix} r\cos(\theta) & -r\sin(\theta) & & & & \\ r\sin(\theta) & r\cos(\theta) & & & & \\ 1 & 0 & r\cos(\theta) & -r\sin(\theta) & & \\ 0 & 1 & r\sin(\theta) & r\cos(\theta) & & \\ & & 1 & 0 & r\cos(\theta) & -r\sin(\theta) \\ & & 0 & 1 & r\sin(\theta) & r\cos(\theta) \end{bmatrix} = \begin{bmatrix} rR(\theta) & 0 & 0 \\ I & rR(\theta) & 0 \\ 0 & I & rR(\theta) \end{bmatrix}$$

Powers of these matrices follow the same pattern as for the $J_k(\lambda)$ blocks (see Proposition 7.1 and the remarks preceding it) but with the observation that $R(\theta)^n = R(n\theta)$. Thus, we have

$$J_3(r,\theta)^n = \begin{bmatrix} r^n R(n\theta) & 0 & 0 \\ nr^{n-1}R((n-1)\theta) & r^n R(n\theta) & 0 \\ \frac{n(n-1)}{2} r^{n-2} R((n-2)\theta) & nr^{n-1}R((n-1)\theta) & r^n R(n\theta) \end{bmatrix}$$

It is known that $C_k(r,\theta)$ on $\mathbb{C}^{2k}$ is similar to $J_k(\lambda) \oplus J_k(\overline{\lambda})$ on $\mathbb{C}^{2k}$ where $\lambda = re^{i\theta}$, see [10, p. 150]. Also, every real $n \times n$ matrix $T$ is similar to its *real Jordan Form* which is a direct sum of blocks of the form $J_k(\lambda)$ where $\lambda$ is a real eigenvalue for $T$ and a direct sum of blocks of the form $C_k(r,\theta)$ where $[r\cos(\theta) \pm ir\sin(\theta)]$ is a conjugate pair of complex eigenvalues for $T$. A *Jordan matrix* is any matrix that is a direct sum of Jordan blocks. For more information on the real Jordan form see [9, p. 359] or [10, p. 150].

**Definition 8.1.** Let $\mathbb{C}^n_{\mathbb{R}}$ denote the set $\mathbb{C}^n$ considered as a vector space over the field $\mathbb{R}$ of real numbers. Then $\mathbb{C}^n_{\mathbb{R}}$ is a $2n$ dimensional (real) vector space. In fact, $\{\vec{e}_k\}_{k=1}^n \cup \{i\vec{e}_k\}_{k=1}^n$ is an orthonormal basis for $\mathbb{C}^n_{\mathbb{R}}$ where $\{\vec{e}_k\}_{k=1}^n$ is the standard unit vector basis for $\mathbb{R}^n$. Also, let $U_c : \mathbb{R}^{2n} \to \mathbb{C}^n_{\mathbb{R}}$ be *the complexification map* given by

$$U_c(x_1, x_2, \ldots, x_{2n-1}, x_{2n}) = (x_1 + ix_2, x_3 + ix_4, \ldots, x_{2n-1} + ix_{2n}).$$

**Proposition 8.2.** ***The Complexification Map & Jordan Blocks***

*If $U_c : \mathbb{R}^{2n} \to \mathbb{C}^n_{\mathbb{R}}$ is the complexification map, then the following hold:*

1. *$U_c$ is a (real) linear isometry mapping $\mathbb{R}^{2n}$ onto $\mathbb{C}^n_{\mathbb{R}}$.*
2. *$U_c C_n(r,\theta) = J_n(\lambda) U_c$ where $\lambda = re^{i\theta}$.*
3. *If $A$ is a $(2n)\times(2n)$ real matrix and $B$ is an $n \times n$ complex matrix and if $U_c A = BU_c$, then $A$ is convex-cyclic on $\mathbb{R}^{2n}$ if and only if $B$ is convex-cyclic on $\mathbb{C}^n_{\mathbb{R}}$ if and only if $B$ is convex-cyclic on $\mathbb{C}^n$.*

*Proof.* Property (1) is elementary. For (2) one may easily verify that $U_c C_n(r,\theta) = J_n(\lambda)U_c$ by checking that $U_c C_n(r,\theta)\vec{e}_k = J_n(\lambda) U_c \vec{e}_k$ for $1 \le k \le 2n$ where $\{\vec{e}_k\}$ is the standard unit vector basis for $\mathbb{R}^{2n}$. For (3), $A$ is convex-cyclic on $\mathbb{R}^{2n}$ if and only if $B$ is convex-cyclic on $\mathbb{C}^n_{\mathbb{R}}$ since $U_c A = BU_c$ holds and convex-cyclicity only involves polynomials with real-coefficients. Lastly, a set $X$ is dense in $\mathbb{C}^n$ if and only if $X$ is dense in $\mathbb{C}^n_{\mathbb{R}}$ since the two *sets* $\mathbb{C}^n$ and $\mathbb{C}^n_{\mathbb{R}}$ are the same and have the same metric, thus the same topologies. Thus the convex-hull (which only involves real scalars) of an orbit produces the same set in both $\mathbb{C}^n$ and $\mathbb{C}^n_{\mathbb{R}}$ and density in $\mathbb{C}^n$ is equivalent to density in $\mathbb{C}^n_{\mathbb{R}}$. □

**Theorem 8.3.** ***Real Jordan Matrices with Complex Eigenvalues.***

1. *If $r \ge 0$, $\theta \in \mathbb{R}$, and $\lambda = re^{i\theta}$, then the real Jordan block $C_n(r,\theta)$ is convex-cyclic on $\mathbb{R}^{2n}$ if and only if the Jordan block $J_n(\lambda)$ is convex-cyclic on $\mathbb{C}^n$ if and only if $\lambda \in \mathbb{C} \setminus (\overline{\mathbb{D}} \cup \mathbb{R})$.*

(2) *If $C = \bigoplus_{k=1}^N C_{n_k}(r_k, \theta_k)$ acts on $\mathbb{R}^{2p}$ where $p = \sum_{k=1}^N n_k$, $r_k \geq 0$, and $\theta_k \in \mathbb{R}$, then $C$ is convex-cyclic on $\mathbb{R}^{2p}$ if and only if $J = \bigoplus_{k=1}^N J_{n_k}(\lambda_k)$ is convex-cyclic on $\mathbb{C}^p$ where $\lambda_k = r_k e^{i\theta_k}$ for $1 \leq k \leq N$ if and only if the following hold:*
   (a) *the eigenvalues $\{\lambda_k\}_{k=1}^N$ are distinct and not real;*
   (b) *$|\lambda_k| > 1$ for all $1 \leq k \leq N$;*
   (c) *for any $1 \leq j, k \leq N$, $\lambda_j \neq \overline{\lambda_k}$.*

   *Furthermore, the convex-cyclic vectors for $C$ are precisely those vectors*

$$\vec{v} = (\vec{v}_1, \vec{v}_2, \ldots, \vec{v}_{p-1}, \vec{v}_p)$$

   *where $\vec{v}_k = (v_{k,1}, v_{k,2}, \ldots, v_{k,2n_k}) \in \mathbb{R}^{2n_k}$ satisfies that for every $1 \leq k \leq p$, $(v_{k,1}, v_{k,2}) \neq (0, 0)$.*

*Proof.* (1) This follows directly from Proposition 7.2 and Proposition 8.2. (2) If $C = \bigoplus_{k=1}^N C_{n_k}(r_k, \theta_k)$ and $J = \bigoplus_{k=1}^N J_{n_k}(r_k e^{i\theta_k})$, then $U_c C = J U_c$ where $U_c$ is the complexification map. Thus by Proposition 8.2, $C$ is convex-cyclic on $\mathbb{R}^{2p}$ if and only if $J$ is convex-cyclic on $\mathbb{C}^p$. The theorem now follows from Theorem 7.5. □

**Theorem 8.4.** ***Real Matrices with Diagonal Complexification.***
*If $D = diag(x_1, x_2, \ldots, x_M)$ is a diagonal matrix on $\mathbb{R}^M$ and $C = \bigoplus_{k=1}^N C_1(r_k, \theta_k)$ acts on $\mathbb{R}^{2N}$, then $T = D \oplus C$ is convex-cyclic on $\mathbb{R}^M \oplus \mathbb{R}^{2N}$ if and only if the following hold, where $\lambda_k = r_k e^{i\theta_k}$ for $1 \leq k \leq N$:*

(1) *the complex eigenvalues $\{\lambda_k\}_{k=1}^N$ are distinct and not real;*
(2) *$|\lambda_k| > 1$ for all $1 \leq k \leq N$;*
(3) *for any $1 \leq j, k \leq N$, $\lambda_j \neq \overline{\lambda_k}$;*
(4) *The $\{x_k\}_{k=1}^M$ are distinct and $x_k < -1$ for all $1 \leq k \leq M$.*

*Furthermore, the convex-cyclic vectors for $T$ are precisely those vectors $\vec{v} = (v_1, v_2, \ldots, v_M, \vec{u}_1, \vec{u}_2, \ldots, \vec{u}_N)$ where $v_j \neq 0$ for all $1 \leq j \leq M$ and $\vec{u}_k \in \mathbb{R}^2 \setminus \{(0, 0)\}$ for all $1 \leq k \leq N$.*

*Proof.* Let $V = I \oplus U_c$ be the mapping $\mathbb{R}^M \oplus \mathbb{R}^{2N} \to \mathbb{R}^M \oplus \mathbb{C}^M_{\mathbb{R}}$. Then $V$ is a real linear onto isometry and $VTV^{-1} = diag(x_1, x_2, \ldots, x_M, \lambda_1, \ldots, \lambda_N)$ is a diagonal matrix. Since (1) - (4) hold, Theorem 4.1 produces peaking convex-polynomials for arbitrary subsets of $\{x_1, \ldots, x_M, \lambda_1, \ldots, \lambda_N\}$ and thus we may apply the same proof from Theorem 5.1 to the case at hand and conclude that $VTV^{-1}$ is convex-cyclic and its convex-cyclic vectors are all vectors with all nonzero coordinates. It follows then that $T$ is convex-cyclic and has the stated set of convex-cyclic vectors. □

**Corollary 8.5.** ***Interpolating Prescribed Values.*** *Suppose that $\{x_k\}_{k=1}^M \subseteq \mathbb{R}$ and $\{z_k\}_{k=1}^N \subseteq \mathbb{C}$ and $z_k = r_k e^{i\theta_k}$ where $r_k \geq 0$ and $\theta_k \in \mathbb{R}$ for all $1 \leq k \leq N$. Suppose also that the following hold:*

(1) *the numbers $\{z_k\}_{k=1}^N$ are distinct and not real;*
(2) *$|z_k| > 1$ for all $1 \leq k \leq N$;*
(3) *for any $1 \leq j, k \leq N$, $z_j \neq \overline{z_k}$.*
(4) *The numbers $\{x_k\}_{k=1}^M$ are distinct and $x_k < -1$ for all $1 \leq k \leq M$.*

*Then given any $\{y_k\}_{k=1}^M \subseteq \mathbb{R}$ and $\{w_k\}_{k=1}^N \subseteq \mathbb{C}$, there exists a convex-polynomial $p$ such that $p(x_k) = y_k$ for $1 \leq k \leq M$ and $p(z_k) = w_k$ for $1 \leq k \leq N$.*

*Proof.* This proof is similar to that of Theorem 6.2, but uses Theorem 8.4 instead of Theorem 5.1 together with the complexification map. □

Notice that in Theorem 8.4 the matrix $C$ was a direct sum of $2\times 2$ blocks. In the theorem below the matrix $C$ is a direct sum of blocks with size $2n_k \times 2n_k$.

**Theorem 8.6.** ***Convex-Cyclicity of $D \oplus C$.***

*If $D = diag(x_1, x_2, \ldots, x_M)$ is a diagonal matrix on $\mathbb{R}^M$ and $C = \bigoplus_{k=1}^N C_{n_k}(r_k, \theta_k)$ on $\mathbb{R}^{2p}$ where $p = \sum_{k=1}^N n_k$ and we let $T = D \oplus C$ on $\mathbb{R}^{M+2p}$, then $T$ is convex-cyclic on $\mathbb{R}^{M+2p}$ if and only if the following hold, where $\lambda_k = r_k e^{i\theta_k}$ for $1 \le k \le N$,*

1. *the complex eigenvalues $\{\lambda_k\}_{k=1}^N$ are distinct and not real;*
2. *$|\lambda_k| > 1$ for all $1 \le k \le N$;*
3. *for any $1 \le j,k \le N$, $\lambda_j \ne \overline{\lambda_k}$.*
4. *The $\{x_k\}_{k=1}^M$ are distinct and $x_k < -1$ for all $1 \le k \le M$.*

*Furthermore, the convex-cyclic vectors for $T$ are precisely those vectors $\vec{v} = (\vec{v}_1, \vec{v}_2)$ where $\vec{v}_1 \in \mathbb{R}^M$ is any convex-cyclic vector for $D$ and $\vec{v}_2 \in \mathbb{R}^{2p}$ is any convex-cyclic vector for $C$.*

*Proof.* We shall use Proposition 3.3 about direct sums of convex-cyclic operators. Given our hypothesis, we know from Theorem 5.1 that $D$ is convex-cyclic and from Theorem 8.3 that $C$ is convex-cyclic. Also, by Corollary 8.5, there exists a convex-polynomial $p$ such that $p(x_k) = -2^k$ for $1 \le k \le M$ and so that $p(\lambda_k) = 0$ for $1 \le k \le N$. Using Proposition 7.1 and property (2) of Proposition 8.2 we see that $p(C)$ is a nilpotent matrix and thus is power bounded. Also, $p(D)$ is a diagonal matrix with diagonal entries $-2^k$ for $1 \le k \le M$, which is convex-cyclic on $\mathbb{R}^M$, by Theorem 5.1. So, Proposition 3.3 implies that $p(T)$ is convex-cyclic and the convex-cyclic vectors for $T$ are direct sums of convex-cyclic vectors as described in the theorem. □

**Corollary 8.7.** ***Interpolating Real Values & Complex Derivatives***

*Suppose that $\{x_k\}_{k=1}^M \subseteq \mathbb{R}$ and $\{z_k\}_{k=1}^N \subseteq \mathbb{C}$ and that the following hold:*

1. *the numbers $\{x_k\}_{k=1}^M$ are distinct and $x_k < -1$ for all $1 \le k \le M$;*
2. *the numbers $\{z_k\}_{k=1}^N$ are distinct and not real;*
3. *$|z_k| > 1$ for all $1 \le k \le N$;*
4. *for any $1 \le j,k \le N$, $z_j \ne \overline{z_k}$.*

*Then given any set $\{y_k\}_{k=1}^M \subseteq \mathbb{R}$ and any set $\{w_{j,k} : 0 \le j \le n, 1 \le k \le N\}$ of complex numbers, there exists a convex-polynomial $p$ such that $p(x_k) = y_k$ for all $1 \le k \le M$ and $p^{(j)}(z_k) = w_{j,k}$ for all $0 \le j \le n$ and $1 \le k \le N$.*

*Proof.* The proof uses Theorem 8.6 and is similar to the proof of Theorem 7.6. □

We are now prepared to show when a real Jordan matrix with real and complex eigenvalues is convex-cyclic.

**Theorem 8.8.** ***Convex-Cyclicity of Real Jordan Matrices.***

*Let $C = \bigoplus_{k=1}^N C_{n_k}(r_k, \theta_k)$ on $\mathbb{R}^{2p}$ where $p = \sum_{k=1}^N n_k$ and $J = \bigoplus_{k=1}^M J_{m_k}(x_k)$ on $\mathbb{R}^q$ where $q = \sum_{k=1}^M m_k$ and let $T = C \oplus J$. Then $T$ is convex-cyclic on $\mathbb{R}^{2p+q}$ if and only if the following hold, where $\lambda_k = r_k e^{i\theta_k}$ for $1 \le k \le N$,*

1. *the complex eigenvalues $\{\lambda_k\}_{k=1}^N$ are distinct and not real;*
2. *$|\lambda_k| > 1$ for all $1 \le k \le N$;*
3. *for any $1 \le j,k \le N$, $\lambda_j \ne \overline{\lambda_k}$.*
4. *The $\{x_k\}_{k=1}^M$ are distinct and $x_k < -1$ for all $1 \le k \le M$.*

*Furthermore, the convex-cyclic vectors for $T$ are precisely those vectors $\vec{v} = (\vec{v}_1, \vec{v}_2)$ where $\vec{v}_1 \in \mathbb{R}^{2p}$ is any convex-cyclic vector for $C$ and $\vec{v}_2 \in \mathbb{R}^q$ is any convex-cyclic vector for $J$.*

*Proof.* We shall use Proposition 3.3 about direct sums of convex-cyclic operators. Given our hypothesis, we know from Theorem 7.5 that $J$ is convex-cyclic and from Theorem 8.3 that $C$ is convex-cyclic. Also by Corollary 8.7, there exists a convex-polynomial $p$ such that $p(x_k) = 0$ for all $1 \leq k \leq M$ and such that $p(\lambda_k) = 2^k i$ and $p'(\lambda_k) = 1$ and $p^{(j)}(\lambda_k) = 0$ for $1 \leq k \leq N$ and $2 \leq j \leq q$. It then follows from Theorem 7.1 that $p(J)$ is nilpotent and thus $p(J)$ is power bounded. Also, $p(C)$ is convex-cyclic. To see this, note that $C$ is similar to the (complex) Jordan matrix $B = \bigoplus_{k=1}^{N} J_k(\lambda_k)$ (via the complexification map) and by Proposition 7.1 and the properties of $p$ we see that $p(B) = \bigoplus_{k=1}^{n} J_k(2^k i)$. By Theorem 7.5 we know that $p(B)$ is convex-cyclic. It follows that $p(C)$ is convex-cyclic. It now follows from Proposition 3.3 that $T = C \oplus J$ is convex-cyclic and has the stated set of convex-cyclic vectors. □

The following interpolation theorem now follows in a similar manner as the previous interpolation theorems did.

**Remark.** It's interesting that the strongest form of the interpolation theorem comes from understanding the convex-cyclicity of real matrices, not just complex ones. This is because convex-cyclic real matrices can have both real and complex eigenvalues. A complex matrix that is convex-cyclic cannot have real eigenvalues.

**Theorem 8.9.** ***Convex-Polynomial Interpolation***
*If $S = \{x_k\}_{k=1}^m \cup \{z_k\}_{k=1}^n \subseteq \mathbb{C}$ where $\{x_k\}_{k=1}^m \subseteq \mathbb{R}$ and $\{z_k\}_{k=1}^n \subseteq \mathbb{C} \setminus \mathbb{R}$, then following are equivalent:*

*(a) for any finite set $\{y_{j,k} : 0 \leq j \leq N, 1 \leq k \leq m\}$ of real numbers and for any finite set $\{w_{j,k} : 0 \leq j \leq N, 1 \leq k \leq n\}$ of complex numbers there exists a convex-polynomial $p$ such that $p^{(j)}(x_k) = y_{j,k}$ for all $0 \leq j \leq N$ and $1 \leq k \leq m$ and $p^{(j)}(z_k) = w_{j,k}$ for all $0 \leq j \leq N$ and $1 \leq k \leq n$.*

*(b) The real numbers $\{x_k\}_{k=1}^m$ are distinct and satisfy $\{x_k\}_{k=1}^m \subseteq (-\infty, -1)$ and the numbers $\{z_k\}_{k=1}^n$ are distinct, $\{z_k\}_{k=1}^n \subseteq \mathbb{C} \setminus (\overline{\mathbb{D}} \cup \mathbb{R})$ and $z_j \neq \overline{z}_k$ for all $1 \leq j, k \leq n$.*

**Theorem 8.10.** ***Convex-Cyclicity of Matrices***
***The Real Case:*** *If $T$ is a real $n \times n$ matrix, then $T$ is convex-cyclic on $\mathbb{R}^n$ if and only if $T$ is cyclic and its real and complex eigenvalues are contained in $\mathbb{C} \setminus (\overline{\mathbb{D}} \cup \mathbb{R}^+)$. If $T$ is convex-cyclic, then the convex-cyclic vectors for $T$ are the same as the cyclic vectors for $T$ and they form a dense set in $\mathbb{R}^n$.*

***The Complex Case:*** *If $T$ is an $n \times n$ matrix, then $T$ is convex-cyclic on $\mathbb{C}^n$ if and only if $T$ is cyclic and its eigenvalues $\{\lambda_k\}_{k=1}^n$ are all contained in $\mathbb{C} \setminus (\overline{\mathbb{D}} \cup \mathbb{R})$ and satisfy $\lambda_j \neq \overline{\lambda}_k$ for all $1 \leq j, k \leq n$. If $T$ is convex-cyclic, then the convex-cyclic vectors for $T$ are the same as the cyclic vectors for $T$ and they form a dense set in $\mathbb{C}^n$.*

Since every matrix is similar to its Jordan Canonical form, the above theorem follows immediately from Theorem 7.5 and Theorem 8.8. Recall that a real or complex matrix is *cyclic if and only if each eigenvalue appears in exactly one Jordan block in its (real or complex) Jordan form*; which means each eigenvalue has geometric multiplicity one (where the geometric multiplicity is the dimension of the eigenspace and in the case of a complex eigenvalue $\lambda$ for a real matrix is defined as the complex dimension of the complex eigenspace corresponding to that eigenvalue.

8.2. **Invariant Convex Sets.** If $T$ is a continuous linear operator on a locally convex space $X$ and $E$ is a subset of $X$, then $E$ is invariant for $T$ if $T(E) \subseteq E$. In this section we determine when the invariant closed convex sets for a matrix are the same as the invariant closed subspaces for the matrix.

This happens exactly when every "part" of the matrix is convex-cyclic. A *part* of an operator $T$ is any operator of the form $T|\mathcal{M}$ where $\mathcal{M}$ is a closed invariant subspace for $T$. So the parts of $T$ are all the operators obtained by restricting $T$ to one of its invariant subspaces.

**Proposition 8.11.** *If $T$ is a continuous linear operator on a locally convex space $X$, then the following are equivalent:*

1. *For every $x \in X$, the closed convex-hull of the orbit of $x$ is a subspace.*
2. *Every closed invariant convex set for $T$ is an invariant subspace for $T$.*
3. *Every cyclic part of $T$ is convex-cyclic and its convex-cyclic vectors are the same as its cyclic vectors. In other words, if $\mathcal{M}$ is a closed invariant subspace for $T$ and $T|\mathcal{M}$ is cyclic, then $T|\mathcal{M}$ is convex-cyclic and the convex-cyclic vectors for $T|\mathcal{M}$ are the same as the cyclic vectors for $T|\mathcal{M}$.*

*Proof.* First notice that a convex set is a subspace if and only if it is closed under scalar multiplication.

(1) $\Rightarrow$ (2). Suppose that (1) holds and let $C$ be a closed invariant convex set for $T$. To show that $C$ is closed under scalar multiplication, let $x \in C$ and let $K = clco(\{T^n x : n \geq 0\}$. Then $K$ is a closed invariant convex set for $T$ and by (1) $K$ is a subspace, hence closed under scalar multiplication, thus $\mathbb{F} \cdot x \subseteq \mathbb{F} \cdot K \subseteq K \subseteq C$. It follows that $C$ is closed under scalar multiplication and hence is a subspace.

(2) $\Rightarrow$ (3). Suppose that (2) holds and let $\mathcal{M}$ be a closed invariant subspace for $T$ and assume that $T|\mathcal{M}$ is cyclic and we will show that $T|\mathcal{M}$ is convex-cyclic and has the appropriate set of convex-cyclic vectors. Let $x \in \mathcal{M}$ be a cyclic vector for $T|\mathcal{M}$. Then the smallest closed invariant subspace for $T$ that contains $x$ is $X$. Now let $C$ be the closed convex-hull of the orbit of $x$ under $T|\mathcal{M}$. Then $C$ is a closed invariant convex set for $T$ and thus by assumption (2) it follows that $C$ is a subspace. Thus $C$ is a closed invariant subspace for $T$ that contains $x$ and hence must be equal to $\mathcal{M}$. Thus $C = \mathcal{M}$ which implies that $x$ is a convex-cyclic vector for $T|\mathcal{M}$.

(3) $\Rightarrow$ (1). Assume that (3) holds and let $x \in X$ and let $C$ be the closed convex-hull of the orbit of $x$ under $T$. We must show that $C$ is a subspace. Let $\mathcal{M}$ be the closed invariant subspace generated by $x$; that is, the closure of the linear span of the orbit of $x$ under $T$. Then $x$ is a cyclic vector for $T|\mathcal{M}$ and thus by property (3) $T|\mathcal{M}$ is convex-cyclic and $x$ is a convex-cyclic vector for $T|\mathcal{M}$. It follows immediately from the definition of a convex-cyclic vector that $C = \mathcal{M}$, hence $C$ is a subspace. □

**Theorem 8.12.** ***Invariant Convex sets for Matrices***

***The Complex Case:*** *A matrix $T$ acting on $\mathbb{C}^n$ has the property that all of its invariant closed convex-sets are invariant subspaces if and only if the eigenvalues $\{\lambda_k\}_{k=1}^n$ of $T$ are contained in $\mathbb{C} \setminus (\overline{\mathbb{D}} \cup \mathbb{R})$ and satisfy $\lambda_j \neq \overline{\lambda}_k$ for all $1 \leq j, k \leq n$.*

***The Real Case:*** *A matrix $T$ acting on $\mathbb{R}^n$ has the property that all of its invariant closed convex-sets are invariant subspaces if and only if its eigenvalues are contained in $\mathbb{C} \setminus (\overline{\mathbb{D}} \cup \mathbb{R}^+)$.*

*Proof.* Since items (2) and (3) in Proposition 8.11 are equivalent it suffices to verify that condition (3) in Proposition 8.11 is equivalent to the condition on the eigenvalues stated in this theorem; and that is exactly what Theorem 8.10 says. □

## 9. Questions

1. Is there a way to explicitly construct interpolating convex-polynomials?
2. If an Abelian semigroup of matrices is convex-cyclic, must it contain a convex-cyclic matrix?

## References

[1] A. Ayadi, H. Marzougui, *Abelian semigroups of matrices on $\mathbb{C}^n$ and hypercyclicity* Proc. Edinb. Math. Soc. (2) **57** (2014), no. 2, 323–338.

[2] F. Bayart, E. Matheron, Dynamics of linear operators. Cambridge Tracts in Mathematics, 179. Cambridge University Press, Cambridge, 2009.

[3] T. Bermúdez, A. Bonilla, & N.S. Feldman, *On Convex-Cyclic Operators*, preprint.

[4] J. B. Conway, A Course in Functional Analysis, 2nd Edition, Springer-Verlag, 1990.

[5] G. Costakis, D. Hadjiloucas, A. Manoussos, *Dynamics of tuples of matrices* Proc. Amer. Math. Soc. **137** (2009), no. 3, 1025–1034.

[6] G. Costakis & I. Parissis, *Dynamics of tuples of matrices in Jordan form* Oper. Matrices **7** (2013), no. 1, 131–157.

[7] L. Elsner, *On matrices leaving invariant a nontrivial convex set*, Linear Algebra Appl. **42** (1982), 103–107.

[8] K. Grosse-Erdmann, A. Peris, Linear chaos. Universitext. Springer, London, 2011.

[9] I. Gohberg, P. Lancaster, & L. Rodman, *Invariant Subspaces of Matrices with Applications*, Canadian Mathematical Society, John Wiley & Sons, 1986.

[10] R.A. Horn & C.R. Johnson, Matrix Analysis, Cambridge University Press, 1999.

[11] F. León-Saavedra, M. P. Romero de la Rosa, *Powers of convex-cyclic operators*, Abstract and Applied Analysis, volume **2014** (2014), Article ID 631894, 3 pages.

[12] H. Rezaei, *On the convex hull generated by orbit of operators*, Linear Algebra and its Applications, **438** (2013), 4190-4203.

Dept. of Mathematics, Washington and Lee University, Lexington VA 24450
*E-mail address*: feldmanN@wlu.edu

Dept. of Mathematics, Bucknell University, Lewisburg PA 17837
*E-mail address*: pmcguire@bucknell.edu